\documentclass[12pt]{article}
\pdfoutput=1
\usepackage{moreverb}
\usepackage[a4paper,margin=2cm]{geometry}

\usepackage{amssymb,amsmath,amsfonts,amsbsy}
\usepackage{algorithmic,algorithm}
\usepackage{graphicx}
\usepackage{multirow}
\usepackage{color}

\begin{document}

\title{A Nested Dissection Approach to Modeling Transport in Nanodevices:
Algorithms and Applications\thanks{Received: 2 Apr., 2013}}
\author{
U. Hetmaniuk\thanks{Department of Applied Mathematics, University of Washington, Seattle, WA. (\texttt{hetmaniu@uw.edu})} ,
Y. Zhao\thanks{Department of Electrical Engineering, University of Washington, Seattle, WA. (\texttt{zhaoyq@uw.edu})}\;
and M. P. Anantram\thanks{Department of Electrical Engineering, University of Washington, Seattle, WA. (\texttt{anant@uw.edu})}
}
\date{}

\maketitle

\begin{abstract}
Modeling nanoscale devices quantum mechanically is a computationally
challenging problem where new methods to solve the underlying equations
are in a dire need.
In this paper, we present an approach to calculate the charge density
in nanoscale devices, within the context of the non equilibrium Green's
function approach.
Our approach exploits recent advances in using an established graph
partitioning approach.
The developed method has the capability to handle open boundary conditions
that are represented by full self energy matrices required for realistic
modeling of nanoscale devices.
Our method to calculate the electron density 
has a reduced complexity compared to
%is up to ten times faster than 
the established recursive Green's function approach.
As an example, we apply our algorithm to a quantum well superlattice
and a carbon nanotube, which are represented by a continuum and tight
binding Hamiltonian respectively, and demonstrate significant
speed up over the recursive method.
\par
{\bf keywords}: nanodevice, numerical simulation, device modeling, quantum transport, superlattice, nanomaterials, graphene, Green's functions, nanotechnology, nanotransistor, tunneling, design
\end{abstract}

\section{Introduction}

With the advent of smaller nanoelectronic devices, where quantum mechanics
is central to the device operation, and new nanomaterials, such as
nanotubes, graphene, and nanowires, quantum mechanical simulations
have become a necessity. The non-equilibrium Green's function (NEGF)
method \cite{Anantram2008aa,2000SuMi...28..253D,datta2002non} has emerged as a powerful modeling approach for
these nanodevices and nanomaterials. The NEGF method is based on the
self-consistent coupling of Schr\"odinger and Poisson equations and
is designed to capture electron scattering effects with phonons.

A typical NEGF-based simulation solves three
Green's function equations,
\begin{align}
\begin{cases}
\mathbf{A}\left(E\right)\mathbf{G}^{r}\left(E\right) 
& = \; \mathbf{I}
\\
\mathbf{A}\left(E\right)\mathbf{G}^{<}\left(E\right) 
& = \; \boldsymbol{\Sigma}^{<}\left(\mathbf{G}^{r}\left(E\right)\right)^{\dagger}
\\
\mathbf{A}\left(E\right)\mathbf{G}^{>}\left(E\right) 
& = \; \boldsymbol{\Sigma}^{>}\left(\mathbf{G}^{r}\left(E\right)\right)^{\dagger}
\end{cases}\label{eq:defineGless1}
\end{align}
 where the sparse matrix $\mathbf{A}$ is defined by
\begin{equation}
\mathbf{A}=E\mathbf{I}-\mathbf{H}-\boldsymbol{\Sigma}^{r}_{L}-\boldsymbol{\Sigma}^{r}_{R}-\boldsymbol{\Sigma}^{r}_{Phonon}\label{eq:defineA}
\end{equation}
 $\mathbf{G}^{r}(E)$ is called the retarded Green's
function, describing local density of states
and the propagation of electrons injected in the device,
and $\left(\mathbf{G}^{r}(E)\right)^{\dagger}$
its Hermitian conjugate.
$\mathbf{G}^{<}(E)$, the lesser Green's
function, represents the electron correlation function for energy
level $E$; the diagonal elements of $\mathbf{G}^{<}(E)$
represent the electron density per unit energy.
$\mathbf{G}^{>}(E)$, the greater Green's
function, represents the hole correlation function
for energy level $E$, which is proportional to the density
of unoccupied states.
$\mathbf{I}$ is the identity matrix and $\mathbf{H}$
the system Hamiltonian. $\boldsymbol{\Sigma}^{r}_{L}$ and $\boldsymbol{\Sigma}^{r}_{R}$
represent the self-energies due to left and right contact coupling
and $\boldsymbol{\Sigma}^{r}_{Phonon}$ corresponds
to the self-energy governing electron-phonon scattering.
The matrix $\boldsymbol{\Sigma}^{<}$
corresponds to the lesser self-energy
and the matrix $\boldsymbol{\Sigma}^{>}$
to the greater self-energy.
The Green functions are then
incorporated in the coupling between the Schr\"odinger and Poisson
equations \textendash{} see \cite{Anantram2008aa} for further details. The self-consistent
solution of the Schr\"odinger and Poisson equations requires to solve
\eqref{eq:defineGless1} many times until consistency is achieved.
It is well appreciated that the computationally intensive
part of this calculation is solving \eqref{eq:defineGless1}
for the diagonal element of $\mathbf{G}^{<}$
(electron density) and $\mathbf{G}^{>}$ at all energies $E$.
The objective of this paper is to present a new algorithm for accelerating
the solution of \eqref{eq:defineGless1}.

The recursive Green's function method (RGF) 
\cite{Haydock1972aa,Haydock1980aa,Haydock1985aa,Lake1997aa,Sols1989aa} 
is an effective method, often used in practice, to compute
$\mathbf{G}^{r}$, $\mathbf{G}^{<}$, and $\mathbf{G}^{>}$.
For elongated devices, this approach remains
the most efficient. Recently, 
the Hierarchical Schur Complement (HSC)
method \cite{Lin2009ab,Lin2011ac}
and the Fast Inverse using Nested Dissection
(FIND) method \cite{Li2008ab,Li2011ab} exploit the nested dissection method
\cite{George1973aa} to exhibit a significant speedup. The key ideas
behind these two algorithms are to partition the whole matrix $\mathbf{A}$
into blocks for an efficient block LU-factorization. This factorization
is then re-used to fill in all diagonal blocks of the Green's
function and some off-diagonal blocks in a specific order. These two
algorithms are more efficient than RGF and have reduced the operation
count down to a multiple of the cost for a block LU-factorization
of a sparse matrix.
Approximate methods, that are efficient in the ballistic limit,
exist also.
For example, the contact block reduction \cite{Mamaluy2003aa,Mamaluy2005aa}
accelerates the computation by using a limited number of modes
to represent the matrices.
The focus of this paper is to develop an exact method that works
in the presence of scattering.
% with self-energy matrices at least block diagonal.
%

For calculating the lesser Green's function $\mathbf{G}^{<}$,
advanced algorithms for an arbitrary sparse matrix are still in their
infancy\footnote{An efficient
algorithm for calculating $\mathbf{G}^{<}$ is also efficient
for computing $\mathbf{G}^{>}$.}. 
The RGF method remains an effective method, especially for
elongated devices. The extension of FIND \cite{Li2011ab} for $\mathbf{G}^{<}$
yields a reduced asymptotic complexity but the constant in front of
the asymptotic term hinders the reduction in runtime. 
Li et al. \cite{Li2012aa} have recently
proposed a modification of FIND for a significant speedup but
their partitioning of the matrix $\mathbf{A}$ requires some pre-processing.
The contribution of this paper is to present an extension of the HSC
method for calculating diagonal blocks for $\mathbf{G}^{<}$ with
partitions from existing graph partitioning libraries (like, for example,
the package METIS \cite{Karypis:1998:FHQ:305219.305248}).

The rest of the paper is organized as follows. Section 2 will review
exact methods and discuss differences between previous algorithms
and the proposed approach. Section 3 will give a mathematical description
of the new algorithm. Finally, Section 4 describes numerical experiments
to highlight its efficiency. The discussion and analysis in this paper
will focus on two-dimensional problems, while three-dimensional problems
will be illustrated in a future publication.

\section{Review of {\it Eexact} Methods}
%\section{Review of existing methods}

\begin{figure}[htbp]
\centering \includegraphics[width=0.7\textwidth]{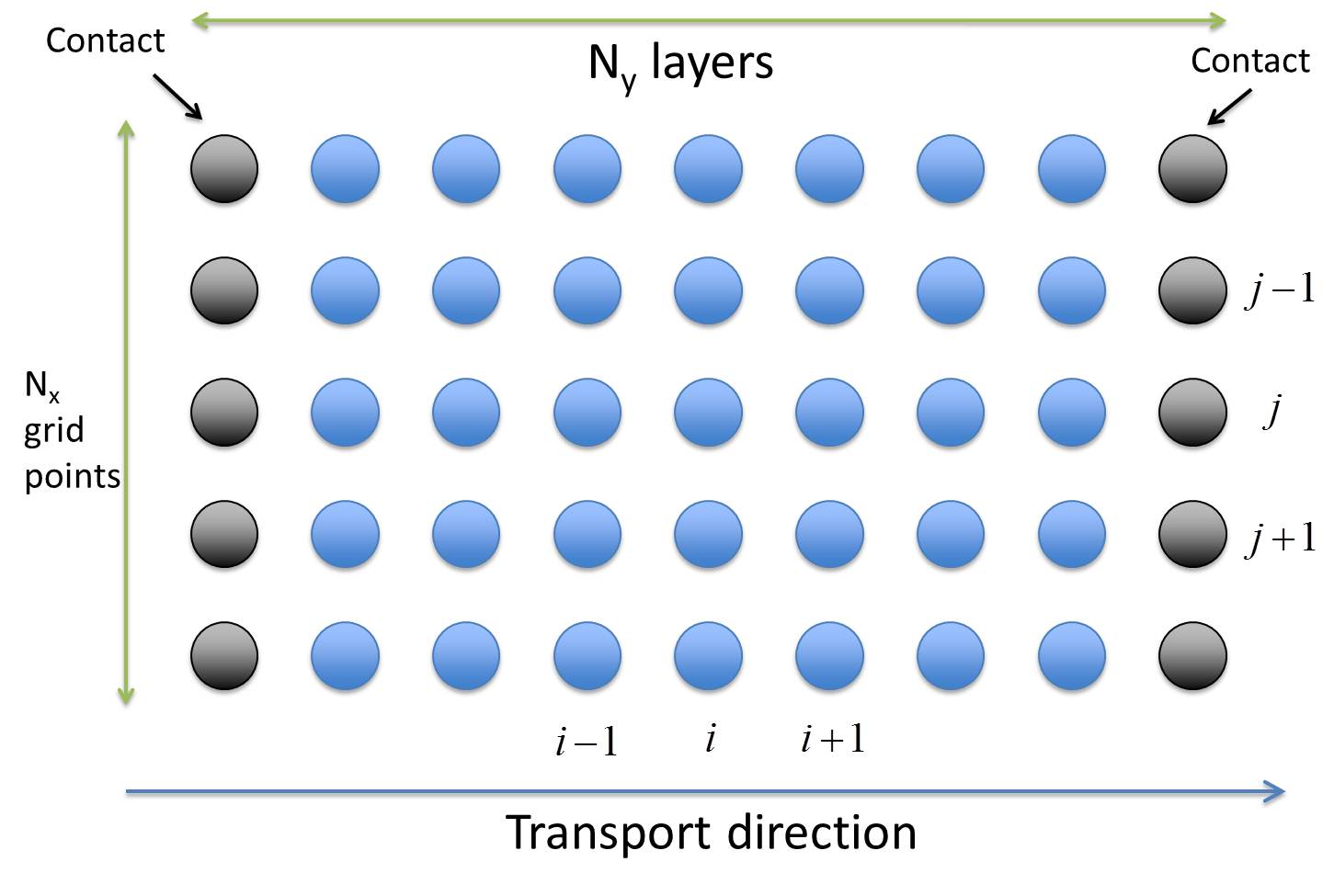}
\caption{Nano-device partitioned into $N_{y}$ layers. Each layer contains
$N_{x}$ grid points.\label{fig:device}}
\end{figure}

Consider a device that can be topologically broken down into layers
as shown in Figure \ref{fig:device}. For an effective mass Hamiltonian,
the blue dots represent grid points of the discretized Green's
function equation, while, in the case of a tight binding Hamiltonian,
the blue dots represent orbitals on an atom.

\begin{figure}[htbp]
\centering \includegraphics[width=0.5\textwidth]{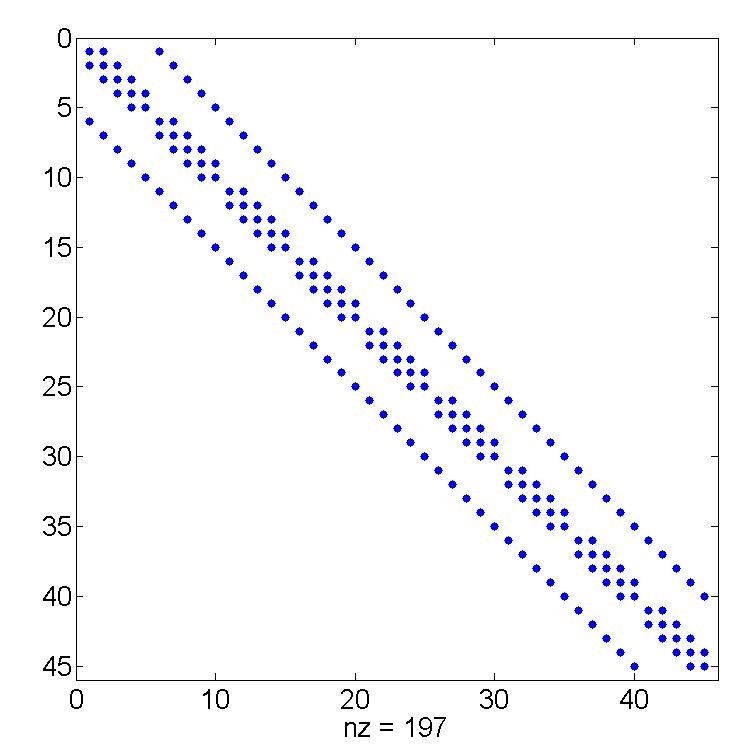}\includegraphics[width=0.5\textwidth]{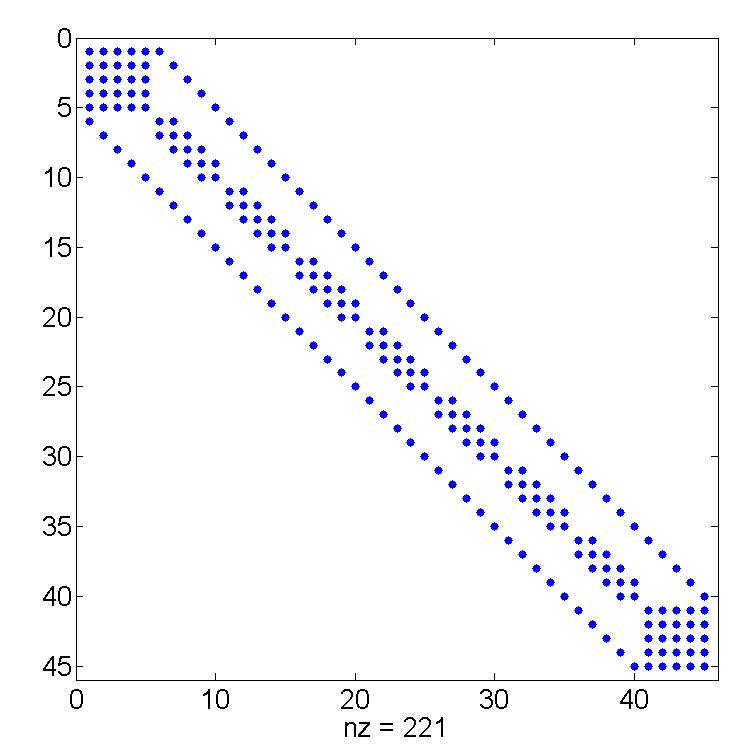}
\caption{$\mathbf{H}$ (left) and $\mathbf{A}$ (right) matrix shape, non-zero
entries are highlighted.\label{fig:sparsity_A_H}}
\end{figure}

When a five-point stencil is used for discretization, the resulting
Hamiltonian is a symmetric block tri-diagonal matrix 
as shown in Figure \ref{fig:sparsity_A_H},
where each diagonal block represents the Hamiltonian of a layer in
Figure \ref{fig:device}. The $i$-th diagonal block of the Hamiltonian
represents the coupling between grid points / atoms in a layer. The
off-diagonal blocks to its left and right represents coupling to between
layers $i$ and $i-1$ and $i+1$ respectively. Both the diagonal
and off-diagonal blocks of the Hamiltonian are sparse in many examples
where either an effective mass or a tight-binding Hamiltonian represents
the dynamics
--- examples of such device include silicon nanowires, nanotubes, and graphene.

The left and right coupling contacts (indicated in Figure \ref{fig:device})
are two semi-infinite leads connected with the device and infinite
matrices represent their Hamiltonians. Their respective effect can
be folded into layer 1 and layer $N_{y}$, resulting in dense blocks
for the first and the $N_{y}$-th diagonal blocks of the self-energy
matrices $\boldsymbol{\Sigma}^{r}_{L}$ and $\boldsymbol{\Sigma}^{r}_{R}$.
The resulting matrix structure of the matrix $\mathbf{A}$, defined
in \eqref{eq:defineA}, is shown in Figure \ref{fig:sparsity_A_H}.
Note that the self-energy matrix $\boldsymbol{\Sigma}^{r}_{Phonon}$
is set to be diagonal at each interior grid point, 
which may arise due to electron-phonon interaction 
or any other interaction.
Relaxing the requirement of diagonal $\boldsymbol{\Sigma}^{r}_{Phonon}$
to include more realistic models
of scattering and solving equations \eqref{eq:defineGless1} 
remains a challenge.

The most common approach to compute blocks of $\mathbf{G}^{r}$ and
$\mathbf{G}^{<}$ is the recursive Green's function
method \cite{Svizhenko2002aa,Anantram2008aa}. RGF is an algorithm
composed of two passes to compute $\mathbf{G}^{r}$ and two passes
to compute $\mathbf{G}^{<}$. In both cases, the passes are interpreted
as follows:
\begin{enumerate}
\item the first pass marches one layer at a time from \textit{left to right}
along the $y$-direction and, recursively, \textit{folds} the effect
of left layers into the current layer;
\item the second pass marches one layer at a time from \textit{right to
left} along the $y$-direction and, recursively, \textit{extracts}
the diagonal blocks and the nearest neighbor off-diagonal blocks for
the final result.
\end{enumerate}
%Further details about RGF are presented in \cite{Svizhenko2002aa,Anantram2008aa}.
Numerically, it is essential to notice that the RGF method exploits
the matrix sparsity of \textit{only} at the block level, which means
that it separates the whole problem into sub-problems of \textit{full}
matrix operations. The complexity of this method is, at most, $10N_{x}^{3}N_{y}$
(when $N_{x}\leq N_{y}$).

To compute block entries of $\mathbf{G}^{r}$, two recent advances,
namely FIND \cite{Li2008ab,Li2011ab} and HSC 
\cite{Lin2009ab,Lin2011ac}, utilize the nested dissection
method \cite{George1973aa} to exhibit a significant speedup. These
methods explicitly exploit the sparsity of $\mathbf{A}$ via a sparse
block LU-factorization of the whole matrix and re-use this factorization
to fill in all diagonal blocks of the Green's function
and some off-diagonal blocks in a specific order. FIND and HSC have
a strong mathematical component and their physical interpretation
is less obvious. The main difference between RGF and these methods
is the replacement of \textit{layers} of grid points organized along
a specific direction with \textit{arbitrarily-shaped clusters} of
grid points organized in a binary tree. Such choice allows to \textit{fold}
and to \textit{extract} in any physical direction when following the
vertical hierarchy of the binary tree. Further details about FIND
and HSC can be found in their respective references. Table \ref{tab:complexity_Gr}
summarizes the complexity of these three state-of-the-art
algorithms when computing entries in $\mathbf{G}^{r}$.

\begin{table}[htbp]
\begin{centering}

\par\end{centering}

\centering{}%
\begin{tabular}{|c|c|c|}
\hline
Algorithm & Complexity when $N_{x}=N_{y}$ & Complexity when $N_{x}<N_{y}$\tabularnewline
\hline
\hline
RGF \cite{Svizhenko2002aa} & $\mathcal{O}\left(N_{x}^{4}\right)$ & $\leq10N_{x}^{3}N_{y}$\tabularnewline
\hline
FIND \cite{Li2008ab,Petersen2009ab} & $\mathcal{O}\left(N_{x}^{3}\right)$ & 
$\leq143N_{x}^{2}N_{y}$
\tabularnewline
\hline
HSC \cite{Lin2009ab,Petersen2009ab} & $\mathcal{O}\left(N_{x}^{3}\right)$ & 
$\leq87N_{x}^{2}N_{y}$
\tabularnewline
\hline
\end{tabular}\caption{Complexity of algorithms to compute diagonal blocks of $\mathbf{G}^{r}$.\label{tab:complexity_Gr}}
\end{table}

Even though FIND \cite{Li2008ab} and HSC \cite{Lin2009ab} have
two distinct mathematical motivations, their runtime complexities
exhibit the same dominating term. But the constants multiplying these
asymptotic terms differ greatly. The runtime for FIND \cite{Li2008ab}
contains a large constant due to the usage of thick boundaries between
the clusters of grid points (or width-2 separators---a width-2 separator
is a boundary between clusters that has a thickness of 2 grid points or 2 atoms).
In addition, generating these clusters with thick boundaries is
not compatible with most existing partitioning libraries. On the other
hand, the HSC method use thin boundaries between clusters with a thickness
of one grid point (or width-1 separators) and can directly exploit
partitions from existing partitioning libraries. As a result, to compute
diagonal blocks for $\mathbf{G}^{r}$, the HSC method \cite[p. 758]{Lin2009ab}
is more efficient than FIND \cite{Li2008ab}.

To compute diagonal entries for $\mathbf{G}^{<}$, only the RGF \cite{Anantram2008aa}
and FIND \cite{Li2011ab} methods have been extended. Table \ref{tab:complexity_Gless}
displays the runtime complexity for computing diagonal blocks of $\mathbf{G}^{<}$.

\begin{table}[htbp]
\centering{}%
\begin{tabular}{|c|c|c|}
\hline
Algorithm & Complexity when $N_{x}=N_{y}$ & Complexity when $N_{x}<N_{y}$\tabularnewline
\hline
\hline
RGF \cite{Svizhenko2002aa} & 
{ $\leq 4N_{x}^{4}$}
 & $\mathcal{O}\left(N_{x}^{3}N_{y}\right)$\tabularnewline
\hline
FIND \cite{Li2011ab} & 
{ $\leq 457 N_{x}^{3}$}
& { $\mathcal{O}\left(N_{x}^{2}N_{y}\right)$}
\tabularnewline
\hline
\end{tabular}\caption{Complexity of algorithms to compute diagonal blocks of $\mathbf{G}^{<}$.\label{tab:complexity_Gless}}
\end{table}

The extension of FIND \cite{Li2011ab} for $\mathbf{G}^{<}$ still
use thick boundaries that result in a large constant for the runtime
and that are incompatible with existing partitioning libraries. 
Recently, Li et al. \cite{Li2012aa} have proposed an extension of
FIND \cite{Li2008ab,Li2011ab} that enables thin boundaries (width-1
separators). This new algorithm yields a significant speedup over
earlier versions of FIND \cite{Li2008ab,Li2011ab}. 
The mathematical motivation remains
different from HSC. The clustering part needs to identify peripheral
sets, an information that is not provided directly by most partitioning
libraries.

The contribution of this paper is to present an extension of the HSC method
\cite{Lin2009ab} for calculating the diagonal blocks for 
$\mathbf{G}^{<}$.
This extension uses thin boundaries (or width-1 separators)
and is compatible with existing partitioning libraries --- namely,
the extension is combined with the graph partitioning package METIS
\cite{Karypis:1998:FHQ:305219.305248}.
The same extension computes efficiently diagonal blocks for 
$\mathbf{G}^{>}$ in a separate step.
For the sake of conciseness, the rest of the paper is focused 
only on computing $\mathbf{G}^{<}$.

\section{Mathematical Description of the Algorithm}
%Possible -> "Proposed method to calculate electron density."

In this section, a detailed mathematical description for the extension
of HSC to compute blocks of $\mathbf{G}^{<}$ is given. The key ingredients
are:
\begin{enumerate}
\item 
an efficient sparse block $\mathbf{LDL}^{T}$-factorization of $\mathbf{A}$.
The block sparse factorization will gather grid points into arbitrarily-shaped
clusters (instead of layers, like in RGF). Such choice allows to fold
and to extract in any physical direction when eliminating entries
in $\mathbf{A}$. The factorization yields formulas to calculate the
diagonal blocks and off-diagonal blocks for $\mathbf{G}^{r}$ and
$\mathbf{G}^{<}$. Exploiting the resulting algebraic relations results
in an algorithm with a cost significantly smaller than the full inversion
of matrix $\mathbf{A}$.
\item an appropriate order of operations. The cost of a matrix multiplication
$\mathbf{BCD}$ depends on the order of operations. When $\mathbf{B}$
is $m\times p$, $\mathbf{C}$ is $p\times k$, and $\mathbf{D}$
is $k\times n$, $\left(\mathbf{BC}\right)\mathbf{D}$ costs $2mk(n+p)$
operations and $\mathbf{B}\left(\mathbf{CD}\right)$ costs $2np(m+k)$.
The order of operations can have a large effect when multiplying series
of matrices together (which is the case for computing entries of $\mathbf{G}^{r}$
and $\mathbf{G}^{<}$). Furthermore, when working with sparse matrices,
one order of operations may preserve sparsity, while another may not.
\end{enumerate}
First, a simple description with three clusters is given. Then the
approach is extended to an arbitrary number of levels and a multilevel
binary tree.

\subsection{Description for a simple case}

\begin{figure}[htbp]
\centering \includegraphics[width=0.7\textwidth]{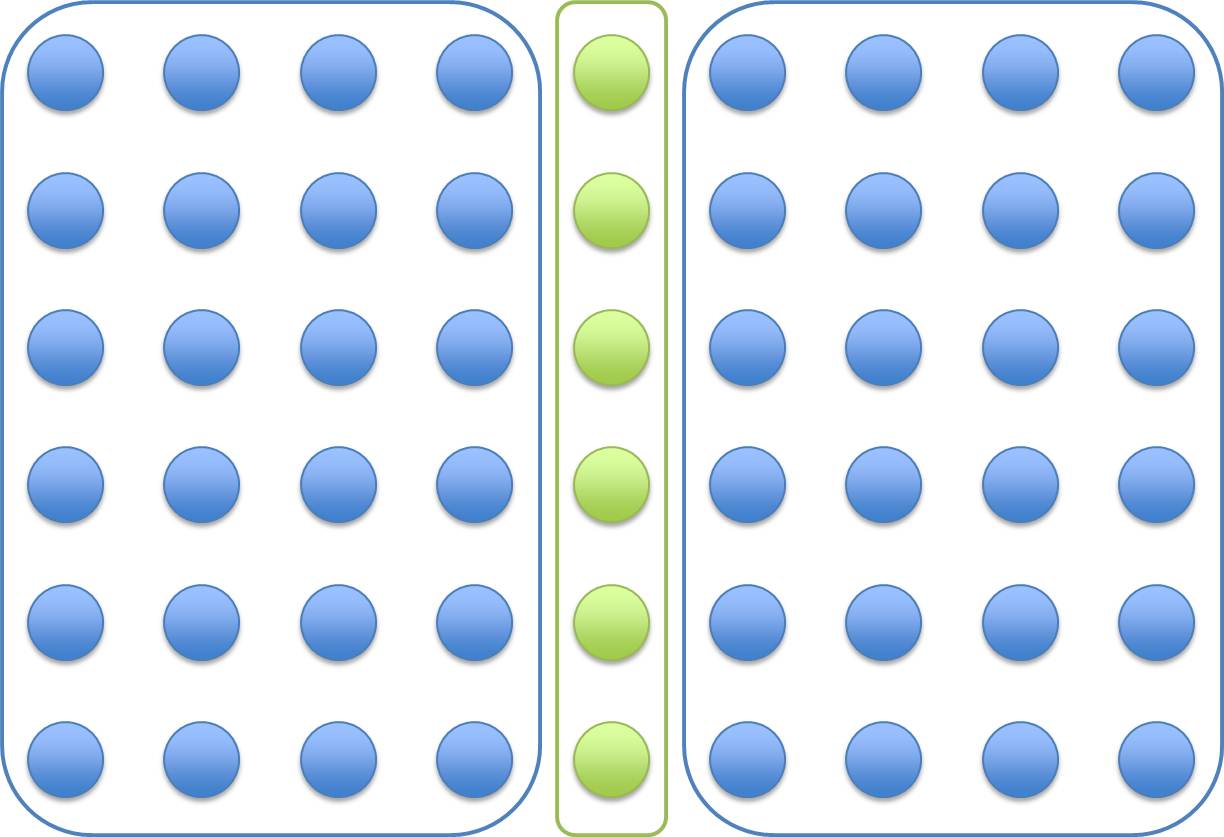}
\caption{Nanodevice partitioned into two subregions (L, R) and a separator
S (L stands for Left and R for Right).\label{fig:device3parts}}
\end{figure}

The basic idea is to partition the nano-device into three disjoint
regions (L, R, S) --- see Figure \ref{fig:device3parts}

Such a partition is easily obtained via the nested dissection, introduced
by George \cite{George1973aa}. Nested dissection divides the system
into two disconnected sets and an interface, called the separator
S. With this partition, the matrix $\mathbf{A}$ can be written as
\[
\mathbf{A}=\left[\begin{array}{ccc}
\mathbf{A}_{LL} & \mathbf{0} & \mathbf{A}_{LS}\\
\mathbf{0} & \mathbf{A}_{RR} & \mathbf{A}_{RS}\\
\mathbf{A}_{LS}^{T} & \mathbf{A}_{RS}^{T} & \mathbf{A}_{SS}
\end{array}\right]
\]
Note that matrix $\mathbf{A}$ is typically complex symmetric.
The block $\mathbf{LDL}^{T}$-factorization of $\mathbf{A}$ is
\[
\mathbf{A}=\left[\begin{array}{ccc}
\mathbf{I} & \mathbf{0} & \mathbf{0}\\
\mathbf{0} & \mathbf{I} & \mathbf{0}\\
\mathbf{A}_{LS}^{T}\mathbf{A}_{LL}^{-1} & \mathbf{A}_{RS}^{T}\mathbf{A}_{RR}^{-1} & \mathbf{I}
\end{array}\right]\left[\begin{array}{ccc}
\mathbf{A}_{LL} & \mathbf{0} & \mathbf{0}\\
\mathbf{0} & \mathbf{A}_{RR} & \mathbf{0}\\
\mathbf{0} & \mathbf{0} & \widehat{\mathbf{A}}_{SS}
\end{array}\right]\left[\begin{array}{ccc}
\mathbf{I} & \mathbf{0} & \mathbf{A}_{LL}^{-1}\mathbf{A}_{LS}\\
\mathbf{0} & \mathbf{I} & \mathbf{A}_{RR}^{-1}\mathbf{A}_{RS}\\
\mathbf{0} & \mathbf{0} & \mathbf{I}
\end{array}\right]
\]
where $\widehat{\mathbf{A}}_{SS}$ is the Schur complement,
$$
\widehat{\mathbf{A}}_{SS}=\mathbf{A}_{SS}-\mathbf{A}_{LS}^{T}\mathbf{A}_{LL}^{-1}\mathbf{A}_{LS}-\mathbf{A}_{RS}^{T}\mathbf{A}_{RR}^{-1}\mathbf{A}_{RS}
.
$$

The matrix $\mathbf{G}^{r}$ satisfies the relation 
\begin{equation}
\mathbf{G}^r
=
\left( \mathbf{I} - \mathbf{L}^{T} \right) \mathbf{G}^r
+ \mathbf{D}^{-1} \mathbf{L}^{-1}
\quad \mbox{with} \quad
\mathbf{A}
=
\mathbf{LD} \mathbf{L}^{T}
\label{eqn:takahashi_gr}
\end{equation}
(described in Takahashi et al. \cite{Takahashi1973} and Erisman and Tinney 
\cite{Erisman1975}).
The block notation yields
\begin{multline*}
\mathbf{G}^{r}=-\left[\begin{array}{ccc}
\mathbf{A}_{LL}^{-1}\mathbf{A}_{LS}\mathbf{G}_{SL}^{r} & \mathbf{A}_{LL}^{-1}\mathbf{A}_{LS}\mathbf{G}_{SR}^{r} & \mathbf{A}_{LL}^{-1}\mathbf{A}_{LS}\mathbf{G}_{SS}^{r}\\
\mathbf{A}_{RR}^{-1}\mathbf{A}_{RS}\mathbf{G}_{SL}^{r} & \mathbf{A}_{RR}^{-1}\mathbf{A}_{RS}\mathbf{G}_{SR}^{r} & \mathbf{A}_{RR}^{-1}\mathbf{A}_{RS}\mathbf{G}_{SS}^{r}\\
\mathbf{0} & \mathbf{0} & \mathbf{\mathbf{0}}
\end{array}\right]\\
+\left[\begin{array}{ccc}
\mathbf{A}_{LL}^{-1} & \mathbf{0} & \mathbf{0}\\
\mathbf{0} & \mathbf{A}_{RR}^{-1} & \mathbf{0}\\
\mathbf{0} & \mathbf{0} & \widehat{\mathbf{A}}_{SS}^{-1}
\end{array}\right]\left[\begin{array}{ccc}
\mathbf{I} & \mathbf{0} & \mathbf{0}\\
\mathbf{0} & \mathbf{I} & \mathbf{0}\\
-\mathbf{A}_{LS}^{T}\mathbf{A}_{LL}^{-1} & -\mathbf{A}_{RS}^{T}\mathbf{A}_{RR}^{-1} & \mathbf{I}
\end{array}\right]
.
\end{multline*}
This equation indicates
\begin{alignat*}{1}
\mathbf{G}_{SS}^{r} & =\left(\widehat{\mathbf{A}}_{SS}\right)^{-1}
,\\
\mathbf{G}_{LS}^{r} & =-\mathbf{A}_{LL}^{-1}\mathbf{A}_{LS}\mathbf{G}_{SS}^{r}
,\\
\mathbf{G}_{RS}^{r} & =-\mathbf{A}_{RR}^{-1}\mathbf{A}_{RS}\mathbf{G}_{SS}^{r}
.
\end{alignat*}
The diagonal blocks $\mathbf{G}_{LL}^{r}$ and $\mathbf{G}_{RR}^{r}$,
for the regions L and R, respectively, are computed independently
of each other,
\begin{alignat*}{1}
\mathbf{G}_{LL}^{r} & =\mathbf{A}_{LL}^{-1}-\mathbf{A}_{LL}^{-1}\mathbf{A}_{LS}\left(\mathbf{G}_{LS}^{r}\right)^{T}=\mathbf{A}_{LL}^{-1}+\mathbf{A}_{LL}^{-1}\mathbf{A}_{LS}\mathbf{G}_{SS}^{r}\mathbf{A}_{LS}^{T}\mathbf{A}_{LL}^{-1}\\
\mathbf{G}_{RR}^{r} & =\mathbf{A}_{RR}^{-1}-\mathbf{A}_{RR}^{-1}\mathbf{A}_{RS}\left(\mathbf{G}_{RS}^{r}\right)^{T}=\mathbf{A}_{RR}^{-1}+\mathbf{A}_{RR}^{-1}\mathbf{A}_{RS}\mathbf{G}_{SS}^{r}\mathbf{A}_{RS}^{T}\mathbf{A}_{RR}^{-1}
\end{alignat*}
(where the symmetry of $\mathbf{G}^{r}$ has been exploited). For
this simple case, the resulting algorithm matches exactly the HSC
method \cite{Lin2009ab}.

For calculating entries in the correlation matrix $\mathbf{G}^{<}$,
Petersen et al. \cite{Petersen2009ab} generalized Takahashi's method
by writing
\begin{equation}
\mathbf{L}^T \mathbf{G}^<  \left( \mathbf{L}^T \right)^\dagger
=
\mathbf{D}^{-1} \mathbf{L}^{-1} \boldsymbol{\Sigma}^{<} 
\mathbf{L}^{-\dagger} \mathbf{D}^{-\dagger} 
\label{eqn:petersen_gless}
. 
\end{equation}
Recurrence formulas are described 
in \cite{Petersen2009ab}.
Our approach utilizes a variant of (\ref{eqn:petersen_gless}), namely
\begin{equation}
\mathbf{G}^{<}
=
\left( \mathbf{I} - \mathbf{L}^{T} \right) \mathbf{G}^{<}
+ \mathbf{D}^{-1} \mathbf{L}^{-1} \boldsymbol{\Sigma}^{<} \left( \mathbf{G}^r \right)^\dagger
\label{eqn:our_gless}
,
\end{equation}
and exploits the symmetry of $\mathbf{G}^{r}$ and
the skew-Hermitian property of $\boldsymbol{\Sigma}^{<}$ and $\mathbf{G}^{<}$,
\begin{equation}
\left( \mathbf{G}^{<} \right)^\dagger = - \mathbf{G}^{<}
\quad \mbox{ and } \quad
\left( \boldsymbol{\Sigma}^{<} \right)^\dagger = - \boldsymbol{\Sigma}^{<}
.
\end{equation}
The block $\mathbf{LDL}^{T}$-factorization of $\mathbf{A}$ yields
\begin{multline*}
\mathbf{G}^{<}=-\left[\begin{array}{ccc}
\mathbf{A}_{LL}^{-1}\mathbf{A}_{LS}\mathbf{G}_{SL}^{<} & \mathbf{A}_{LL}^{-1}\mathbf{A}_{LS}\mathbf{G}_{SR}^{<} & \mathbf{A}_{LL}^{-1}\mathbf{A}_{LS}\mathbf{G}_{SS}^{<}\\
\mathbf{A}_{RR}^{-1}\mathbf{A}_{RS}\mathbf{G}_{SL}^{<} & \mathbf{A}_{RR}^{-1}\mathbf{A}_{RS}\mathbf{G}_{SR}^{<} & \mathbf{A}_{RR}^{-1}\mathbf{A}_{RS}\mathbf{G}_{SS}^{<}\\
\mathbf{0} & \mathbf{0} & \mathbf{\mathbf{0}}
\end{array}\right]\\
+\left[\begin{array}{ccc}
\mathbf{A}_{LL}^{-1} & \mathbf{0} & \mathbf{0}\\
\mathbf{0} & \mathbf{A}_{RR}^{-1} & \mathbf{0}\\
\mathbf{0} & \mathbf{0} & \left(\widehat{\mathbf{A}}_{SS}\right)^{-1}
\end{array}\right]\left[\begin{array}{ccc}
\mathbf{I} & \mathbf{0} & \mathbf{0}\\
\mathbf{0} & \mathbf{I} & \mathbf{0}\\
-\mathbf{A}_{LS}^{T}\mathbf{A}_{LL}^{-1} & -\mathbf{A}_{RS}^{T}\mathbf{A}_{RR}^{-1} & \mathbf{I}
\end{array}\right]\boldsymbol{\Sigma}^{<}
\left( \mathbf{G}^{r} \right)^{\dagger}
.
\end{multline*}
Parts of $\mathbf{G}^{r}$ are computed with the previous algorithm, namely, $\mathbf{G}_{LL}^{r}$,
$\mathbf{G}_{RR}^{r}$, $\mathbf{G}_{SS}^{r}$, $\mathbf{G}_{RS}^{r}$,
and $\mathbf{G}_{LS}^{r}$. 
By assumption, $\boldsymbol{\Sigma}^{<}$ is 
a block-diagonal skew-Hermitian matrix with purely
imaginary entries.
The partial matrix multiplication gives
\begin{multline*}
\left[\begin{array}{ccc}
\mathbf{I} & \mathbf{0} & \mathbf{0}\\
\mathbf{0} & \mathbf{I} & \mathbf{0}\\
-\mathbf{A}_{LS}^{T}\mathbf{A}_{LL}^{-1} & -\mathbf{A}_{RS}^{T}\mathbf{A}_{RR}^{-1} & \mathbf{I}
\end{array}\right]\boldsymbol{\Sigma}^{<} \left( \mathbf{G}^{r} \right)^\dagger
=
\\
\left[\begin{array}{ccc}
\boldsymbol{\Sigma}_{LL}^{<} \left( \mathbf{G}_{LL}^{r} \right)^\dagger 
& \mathbf{*} & \boldsymbol{\Sigma}_{LL}^{<} \left( \mathbf{G}_{SL}^{r} \right)^\dagger 
\\
* & \boldsymbol{\Sigma}_{RR}^{<} \left( \mathbf{G}_{RR}^{r} \right)^\dagger 
& \boldsymbol{\Sigma}_{RR}^{<} \left( \mathbf{G}_{RS}^{r} \right)^\dagger  \\
\mathbf{*} & \mathbf{*} & \boldsymbol{\Sigma}_{SS}^{<} \left( \mathbf{G}_{SS}^{r} \right)^\dagger
-\mathbf{A}_{LS}^{T}\mathbf{A}_{LL}^{-1}\boldsymbol{\Sigma}_{LL}^{<} \left( \mathbf{G}_{SL}^{r} \right)^\dagger
-\mathbf{A}_{RS}^{T}\mathbf{A}_{RR}^{-1}\boldsymbol{\Sigma}_{RR}^{<} \left( \mathbf{G}_{SR}^{r} \right)^\dagger
\end{array}\right]
\end{multline*}
(where the starred blocks are not computed). 
This relation indicates
\begin{alignat*}{1}
\mathbf{G}_{SS}^{<} & =\mathbf{G}_{SS}^{r}\left(
\boldsymbol{\Sigma}_{SS}^{<}\left( \mathbf{G}_{SS}^{r} \right)^\dagger
-\mathbf{A}_{LS}^{T}\mathbf{A}_{LL}^{-1}\boldsymbol{\Sigma}_{LL}^{<}
\left( \mathbf{G}_{SL}^{r} \right)^\dagger 
-\mathbf{A}_{RS}^{T}\mathbf{A}_{RR}^{-1}\boldsymbol{\Sigma}_{RR}^{<}
\left( \mathbf{G}_{SR}^{r} \right)^\dagger
\right)
,\\
\mathbf{G}_{LS}^{<} & =-\mathbf{A}_{LL}^{-1}\mathbf{A}_{LS}\mathbf{G}_{SS}^{<}+\mathbf{A}_{LL}^{-1}\boldsymbol{\Sigma}_{LL}^{<}
\left( \mathbf{G}_{SL}^{r} \right)^\dagger
,\\
\mathbf{G}_{RS}^{<} & =-\mathbf{A}_{RR}^{-1}\mathbf{A}_{RS}\mathbf{G}_{SS}^{<}+\mathbf{A}_{RR}^{-1}\boldsymbol{\Sigma}_{RR}^{<}
\left( \mathbf{G}_{SR}^{r} \right)^\dagger
.
\end{alignat*}
Finally, the diagonal blocks $\mathbf{G}_{LL}^{<}$ and $\mathbf{G}_{RR}^{<}$,
for the regions $L$ and $R$, respectively, are computed independently
of each other,
\begin{alignat*}{1}
\mathbf{G}_{LL}^{<} & =\mathbf{A}_{LL}^{-1}\boldsymbol{\Sigma}_{LL}^{<}
\left( \mathbf{G}_{LL}^{r} \right)^\dagger
-\mathbf{A}_{LL}^{-1}\mathbf{A}_{LS}\left(\mathbf{G}_{LS}^{<}\right)^{\dagger}
\\
\mathbf{G}_{RR}^{<} & =\mathbf{A}_{RR}^{-1}\boldsymbol{\Sigma}_{RR}^{<}
\left( \mathbf{G}_{RR}^{r} \right)^\dagger
-\mathbf{A}_{RR}^{-1}\mathbf{A}_{RS}\left(\mathbf{G}_{RS}^{<}\right)^{\dagger}
\end{alignat*}
(where the skew-Hermitian property of $\mathbf{G}^<$ has been exploited).

%This derivation is the first extension of the HSC method for calculating
%the correlation matrix $\mathbf{G}^{<}$. It differs from the FIND
%method \cite{Li2011ab,Li2012aa} in several aspects. 
This derivation for calculating the correlation matrix $\mathbf{G}^{<}$ 
differs from the FIND method \cite{Li2011ab,Li2012aa} in several aspects. 
It uses thin boundaries obtained directly from the nested dissection. 
This extension requires one sparse factorization and one back-substitution,
while FIND utilizes only sparse factorizations but applied many times
with different orderings. 
The order of operations to obtain the diagonal blocks of $\mathbf{G}^{<}$ 
is also different from the recurrence in Petersen et al. \cite{Petersen2009ab},
which uses the sequence
\begin{multline*}
\boldsymbol{\Sigma}^{<}
\rightarrow
\mathbf{L}^{-1}\boldsymbol{\Sigma}^{<}\left(\mathbf{L}^{-1}\right)^{\dagger}
\rightarrow
\mathbf{D}^{-1}\mathbf{L}^{-1}\boldsymbol{\Sigma}^{<}
\left(\mathbf{L}^{-1}\right)^{\dagger}\left(\mathbf{D}^{-1}\right)^{\dagger}
\\
\rightarrow
\mathbf{L}^{-T}\mathbf{D}^{-1}\mathbf{L}^{-1}\boldsymbol{\Sigma}^{<}
\left(\mathbf{L}^{-1}\right)^{\dagger}\left(\mathbf{D}^{-1}\right)^{\dagger}
\left(\mathbf{L}^{-T}\right)^{\dagger}
,
\end{multline*}
while our HSC extension uses
\[
\boldsymbol{\Sigma}^{<}\rightarrow\boldsymbol{\Sigma}^{<}\left(\mathbf{G}^{r}\right)^{\dagger}\rightarrow\mathbf{L}^{-1}\boldsymbol{\Sigma}^{<}\left(\mathbf{G}^{r}\right)^{\dagger}\rightarrow\mathbf{D}^{-1}\mathbf{L}^{-1}\boldsymbol{\Sigma}^{<}\left(\mathbf{G}^{r}\right)^{\dagger}\rightarrow\mathbf{L}^{-T}\mathbf{D}^{-1}\mathbf{L}^{-1}\boldsymbol{\Sigma}^{<}\left(\mathbf{G}^{r}\right)^{\dagger}
\]
When working with sparse matrices, specific order of operations may
result in fewer operations. 
In numerical experiments, the latter ordering was more efficient.

\subsection{Description for a multilevel case\label{sub:description_multi}}

In this section, the description is extended to an arbitrary number
of clusters.

Even though computing the diagonal for the inverse of a matrix is
not equivalent to a sparse factorization, both problems benefit from
matrix reordering. The multilevel nested dissection, introduced by
George \cite{George1973aa}, lends itself naturally to the creation
of grid points clusters. Typically, nested dissection divides the
system into two disconnected sets and an interface, called the separator.
Then the process is repeated recursively on each set to create a multilevel
binary tree.

\begin{figure}[htbp]
\centering \includegraphics[width=0.7\textwidth]{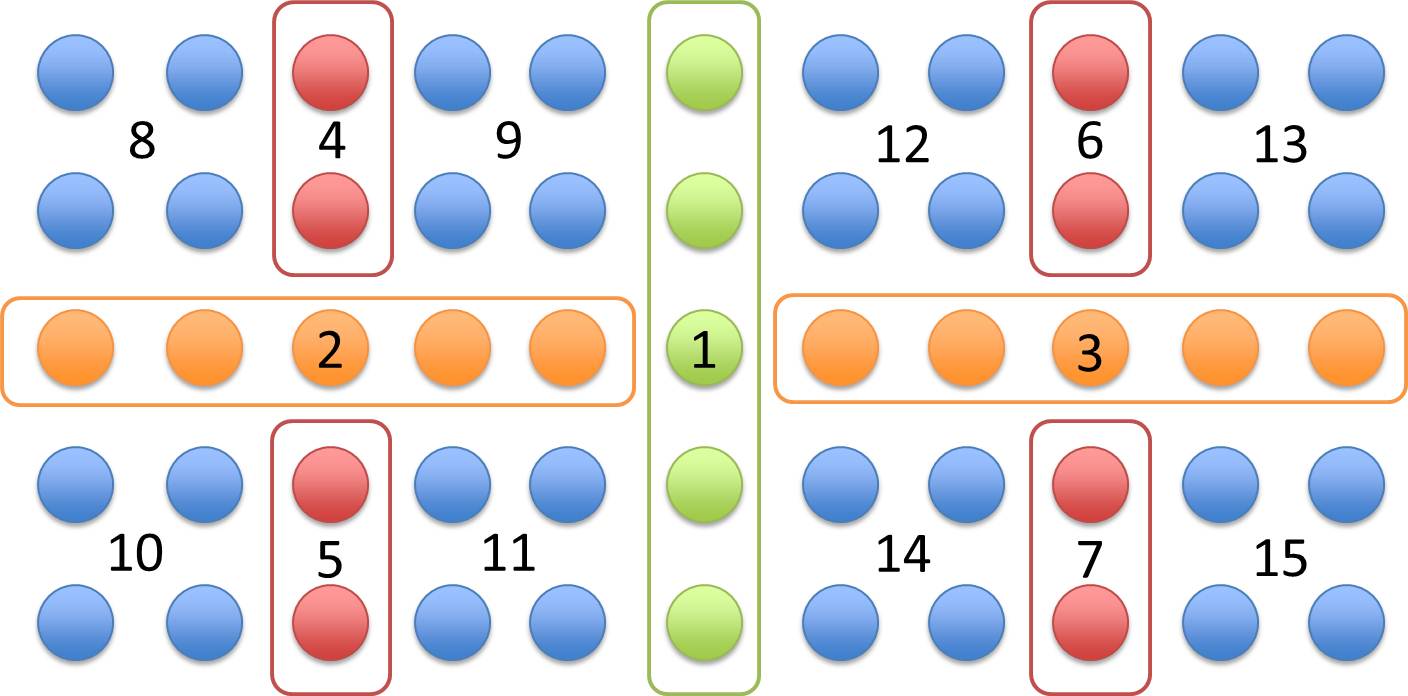}
\centering \includegraphics[width=0.7\textwidth]{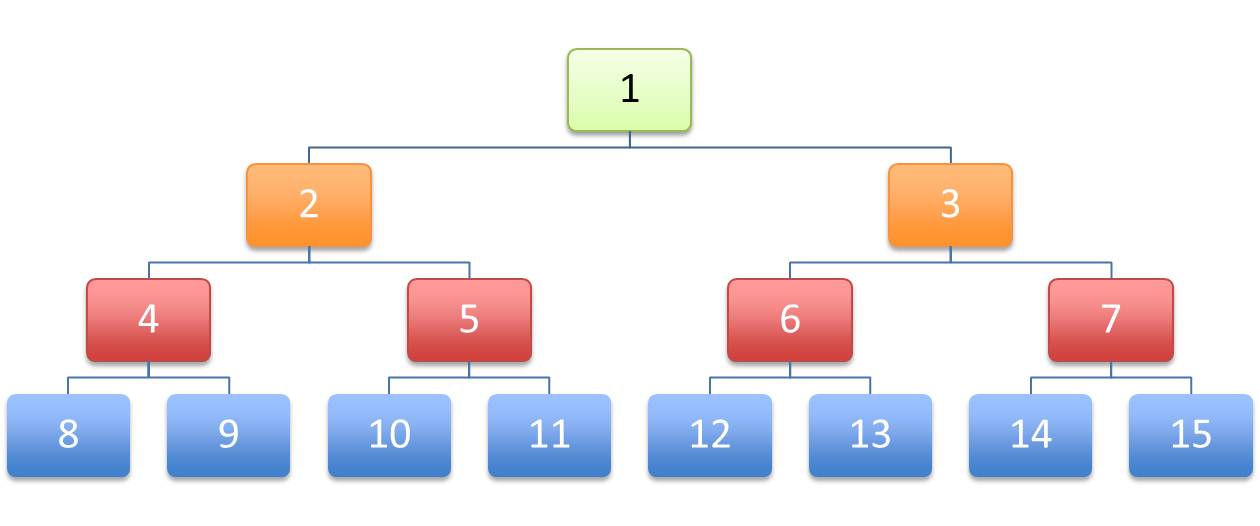}
\caption{Example of a multilevel partition.\label{fig:multilevel_device}}
\end{figure}
Based on the hierarchical structure of the tree, let $P_{i}$ denote
the set of all cluster indices $j$ such that cluster $j$ is an ancestor
of cluster $i$. For example for Figure \ref{fig:multilevel_device},
$P_{5}$ is equal to $\left\{ 1,2\right\} $ and $P_{15}=\left\{ 1,3,7\right\} $.
Let $C_{i}$ denote the set of all cluster indices $j$ such that
cluster $j$ is a descendant of cluster $i$. For the partition on
Figure \ref{fig:multilevel_device}, $C_{4}$ is equal to $\left\{ 8,9\right\} $
and $C_{3}=\left\{ 6,7,12,13,14,15\right\} $. Note that a cluster
may or may not have a direct coupling in the matrix $\mathbf{A}$
to any of its ancestors or descendants.

Once the partition is set, the algorithm may be separated into two
distinct parts: computation of $\mathbf{G}^{r}$ and computation of
$\mathbf{G}^{<}$.

\subsubsection*{Computation of blocks for $\mathbf{G}^{r}$}

In the binary tree, the levels are labeled from bottom to top, where
level 1 contains all clusters at the end of the tree and level $L$
contains only the original separator. For simplicity of presentation,
let $\mathbf{A}^{(l)}$ denote the matrix transformed from $\mathbf{A}$
after folding all the clusters up to level $l$. Note that $\mathbf{A}^{(0)}$
is set to $\mathbf{A}$ and $\mathbf{A}^{(L-1)}$ is block diagonal.

The computation of blocks for $\mathbf{G}^{r}$ involve three steps:
folding the lower level clusters unto the higher ones, inversion of
the matrix for the main separator, and extracting of the diagonal
blocks for the current level from blocks on higher level.

The algorithm for the first step goes as follows:
\begin{itemize}
\item For $l=1$ up to $L-1$,

\begin{itemize}
\item $\mathbf{A}^{(l)}=\mathbf{A}^{(l-1)}$
\item For all the clusters $i$ on level $l$,

\begin{itemize}
\item $\boldsymbol{\Psi}_{i,j}=-\left(\mathbf{A}_{i,i}^{(l)}\right)^{-1}\mathbf{A}_{i,j}^{(l)}$
for all $j$ in $P_{i}$
\item $\mathbf{A}_{j,k}^{(l)}=\mathbf{A}_{j,k}^{(l)}+\boldsymbol{\Psi}_{i,j}^{T}\mathbf{A}_{i,k}^{(l)}$
for all $j$ and $k$ in $P_{i}$
\item $\mathbf{A}_{k,j}^{(l)}=\left(\mathbf{A}_{j,k}^{(l)}\right)^{T}$
for all $j$ and $k$ in $P_{i}$
\item $\mathbf{A}_{i,j}^{(l)}=\mathbf{0}$ and $\mathbf{A}_{j,i}^{(l)}=\mathbf{0}$
for all $j$ in $P_{i}$
\end{itemize}
\item end
\end{itemize}
\item end
\end{itemize}
The next step is written as the inversion of $\mathbf{A}^{(L-1)}$
, which is symmetric and block diagonal.
\begin{itemize}
\item $\mathbf{G}^{(L-1)}=\left(\mathbf{A}^{(L-1)}\right)^{-1}$
\end{itemize}
In practice, the operation requires only the inversion of the block
for the top separator. All the other blocks have already been inverted
during the folding steps.

Finally, all the diagonal blocks of $\mathbf{G}^{r}$ are extracted
one level at a time. The algorithm goes as follows:
\begin{itemize}
\item For $l=L-2$ down to $0$,

\begin{itemize}
\item $\mathbf{G}^{(l)}=\mathbf{G}^{(l+1)}$
\item For all the clusters $i$ on level $l$,

\begin{itemize}
\item $\mathbf{G}_{i,j}^{(l)}=\mathbf{G}_{i,j}^{(l)}+\sum_{k\in P_{i}}\boldsymbol{\Psi}_{i,k}\mathbf{G}_{k,j}^{(l)}$
for all cluster indices $j$ in $P_{i}$
\item $\mathbf{G}_{j,i}^{(l)}=\left(\mathbf{G}_{i,j}^{(l)}\right)^{T}$
for all cluster indices $j$ in $P_{i}$
\item $\mathbf{G}_{i,i}^{(l)}=\mathbf{G}_{i,i}^{(l)}+\sum_{j\in P_{i}}\boldsymbol{\Psi}_{i,j}\mathbf{G}_{j,i}^{(l)}$
\end{itemize}
\item end
\end{itemize}
\item end
\end{itemize}
The resulting algorithm to compute block entries in $\mathbf{G}^{r}$
matches exactly the HSC method \cite{Lin2009ab}.

Note that matrix $\mathbf{G}^{(0)}$ is not equal to $\mathbf{G}^{r}$
because $\mathbf{G}^{(0)}$ is incomplete (see an example in the Appendix).
However, all the entries in $\mathbf{G}^{(0)}$, in particular the
diagonal entries, match the corresponding entries in $\mathbf{G}^{r}$.

\subsubsection*{Computation of blocks for $\mathbf{G}^{<}$}

The algorithm consists of four steps. The first step uses the matrix
$\mathbf{G}^{(0)}$ computed previously.
\begin{itemize}
\item $\mathbf{N}=\boldsymbol{\Sigma}^{<}\left(\mathbf{G}^{(0)}\right)^{\dagger}$
\end{itemize}
All the entries in $\mathbf{G}^{(0)}$ match the corresponding entries
in $\mathbf{G}^{r}$ and are sufficient to compute diagonal blocks
of $\mathbf{G}^{<}$. The matrix $\mathbf{\Sigma}^{<}$ is typically
block diagonal. The matrix $\mathbf{N}$ will have the same structure
and shape as $\mathbf{G}^{(0)}$. The matrix multiplication is done
block by block.

Next the lower level clusters are folded into the higher ones. This
step is critical and the most time consuming. Let $\mathbf{N}^{(l)}$
denote the matrix transformed from $\mathbf{N}$ after folding all
the clusters up to level $l$. $\mathbf{N}^{(0)}$ is set to $\mathbf{N}$.
\begin{itemize}
\item For $l=1$ up to $L-1$,

\begin{itemize}
\item $\mathbf{N}^{(l)}=\mathbf{N}^{(l-1)}$
\item For all the clusters $i$ on level $l$,

\begin{itemize}
\item $\mathbf{N}_{j,k}^{(l)}=\mathbf{N}_{j,k}^{(l)}+\boldsymbol{\Psi}_{i,j}^{T}\mathbf{N}_{i,k}^{(l)}$
for all $j$ and $k$ in $P_{i}$
\end{itemize}
\item end
\end{itemize}
\item end
\end{itemize}
Similarly to Step 1, Step 3 is a block diagonal multiplication.
\begin{itemize}
\item $\mathbf{P}^{(L-1)}=\mathbf{G}^{(L-1)}\mathbf{N}^{(L-1)}$
\end{itemize}
Finally, Step 4 extracts all the diagonal blocks one level at a time.
This step is similar to the extraction in the $\mathbf{G}^{r}$ algorithm.
The operations are the following:
\begin{itemize}
\item For $l=L-2$ down to $0$,

\begin{itemize}
\item $\mathbf{P}^{(l)}=\mathbf{P}^{(l+1)}$
\item For all the clusters $i$ on level $l$,

\begin{itemize}
\item $\mathbf{P}_{i,j}^{(l)}=\mathbf{P}_{i,j}^{(l)}+\sum_{k\in P_{i}}\boldsymbol{\Psi}_{i,k}\mathbf{P}_{k,j}^{(l)}$
for all cluster indices $j$ in $P_{i}$
\item $\mathbf{P}_{j,i}^{(l)}=-\left(\mathbf{P}_{i,j}^{(l)}\right)^{\dagger}$
for all cluster indices $j$ in $P_{i}$
\item $\mathbf{P}_{i,i}^{(l)}=\mathbf{P}_{i,i}^{(l)}+\sum_{j\in P_{i}}\boldsymbol{\Psi}_{i,j}\mathbf{P}_{j,i}^{(l)}$
\end{itemize}
\item end
\end{itemize}
\item end
\end{itemize}
At the end, matrix $\mathbf{P}^{(0)}$ is not equal to $\mathbf{G}^{<}$
because $\mathbf{P}^{(0)}$ is incomplete (see an example in the Appendix).
However, all the entries in $\mathbf{P}^{(0)}$, in particular the
diagonal entries, match the corresponding entries in $\mathbf{G}^{<}$.

\subsubsection*{Comments on the system partition}

The partitioning of the system (or the clustering of points) is the
key step for the efficiency of this algorithm. The partition should
follow two rules:
\begin{enumerate}
\item for clusters within the same level on the binary tree, no interaction
is allowed. Operations on blocks at the same level are performed independently.
\item the partition should minimize the size of separators and reduce the
clusters down to a size manageable for an inversion of the corresponding
block matrix.
\end{enumerate}
The multilevel nested dissection generates a partition that satisfies
those rules. It is worthwhile to note that, as long as the rules stated
above are followed, systems with non-uniform distribution of points
or with a different stencil could be treated correctly.

\begin{figure}[htbp]
\centering \includegraphics[width=0.7\textwidth]{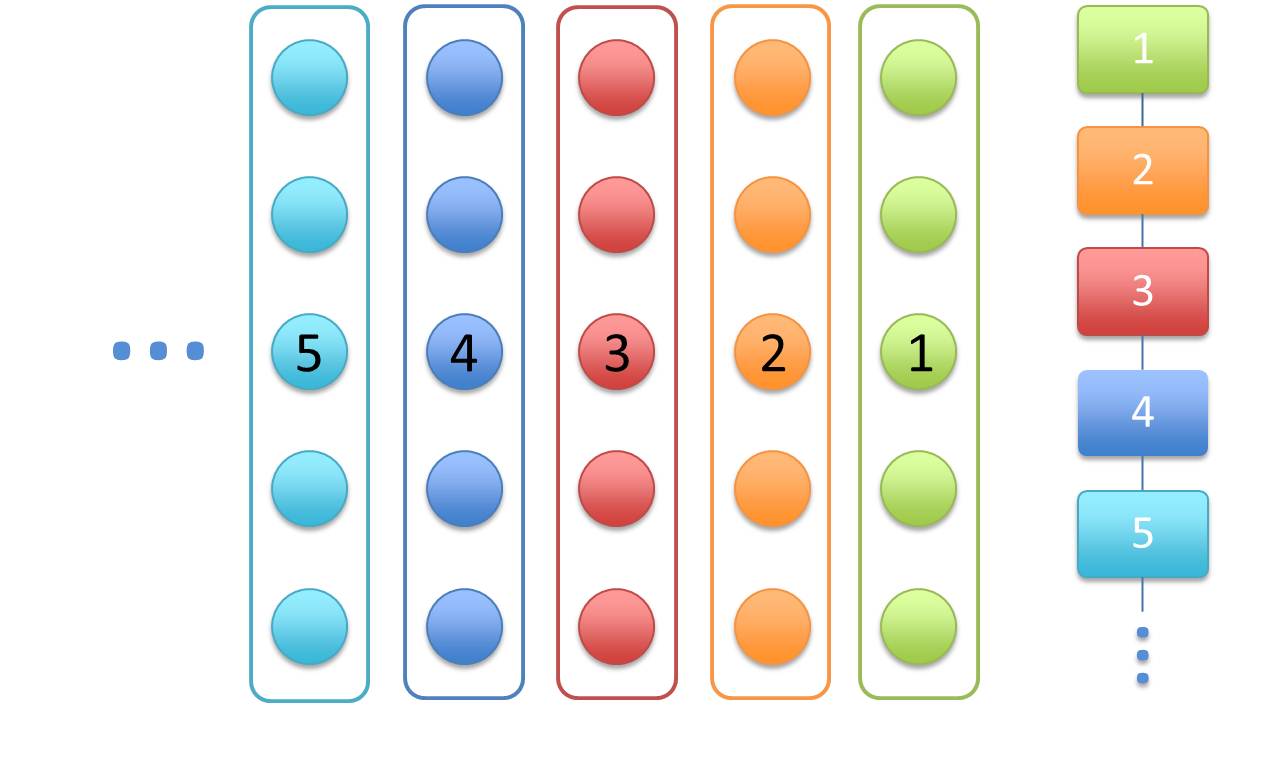}
\caption{Partition generating the RGF algorithm.\label{fig:partition_RGF}}
\end{figure}

The RGF algorithm is included in the previous description with a very
specific partition.
Such a partition is illustrated in Figure \ref{fig:partition_RGF}.

\begin{figure}[htbp]
\centering \includegraphics[width=0.7\textwidth]{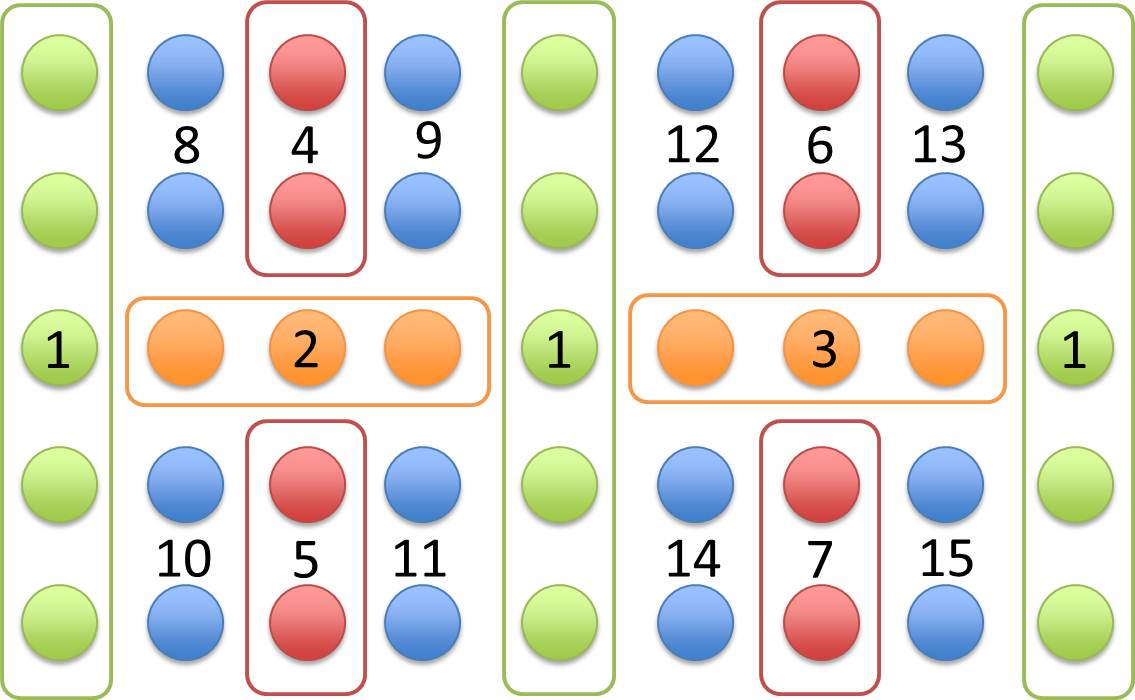}
\caption{First method to partition system with two dense layers and two ends.\label{fig:partition_1}}
\end{figure}

\begin{figure}[htbp]
\centering \includegraphics[width=0.7\textwidth]{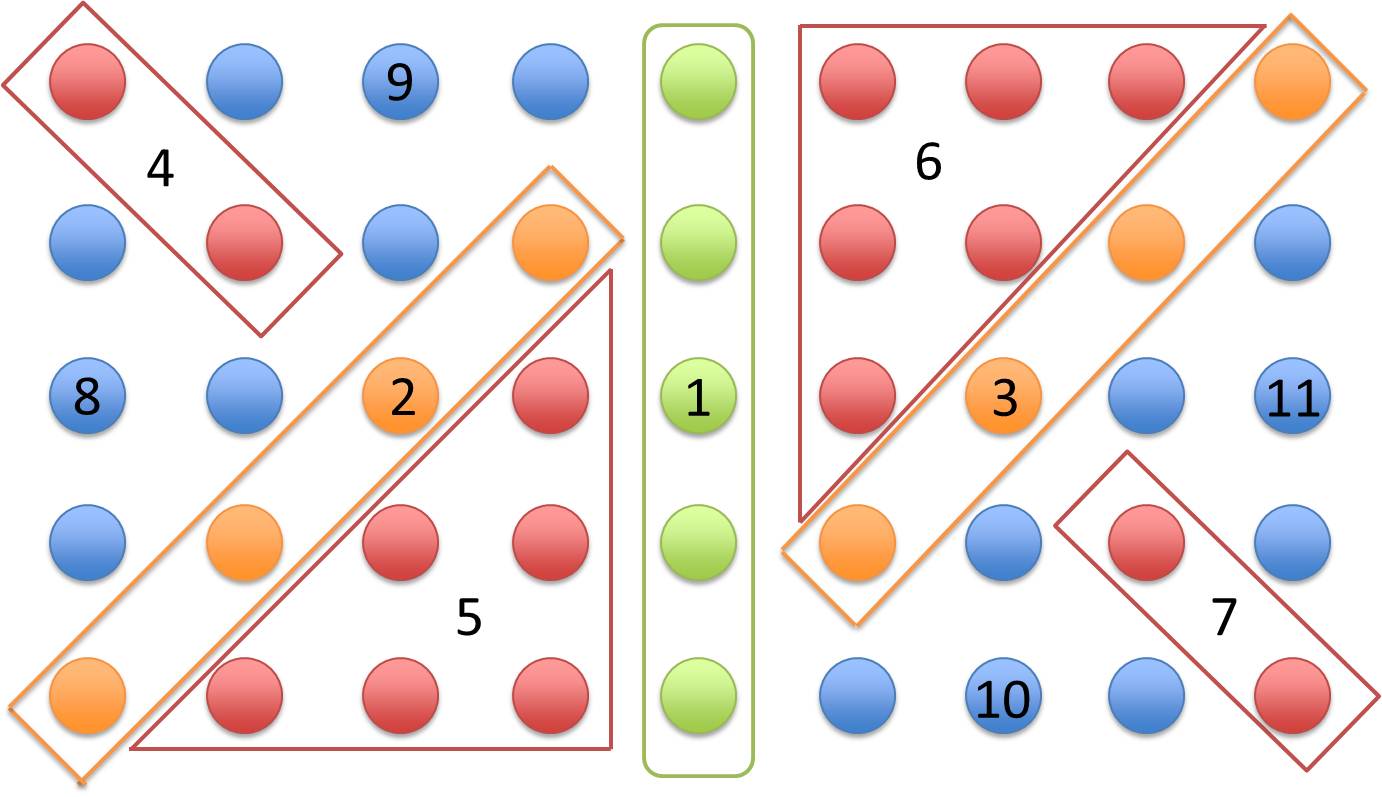}
\includegraphics[width=0.7\textwidth]{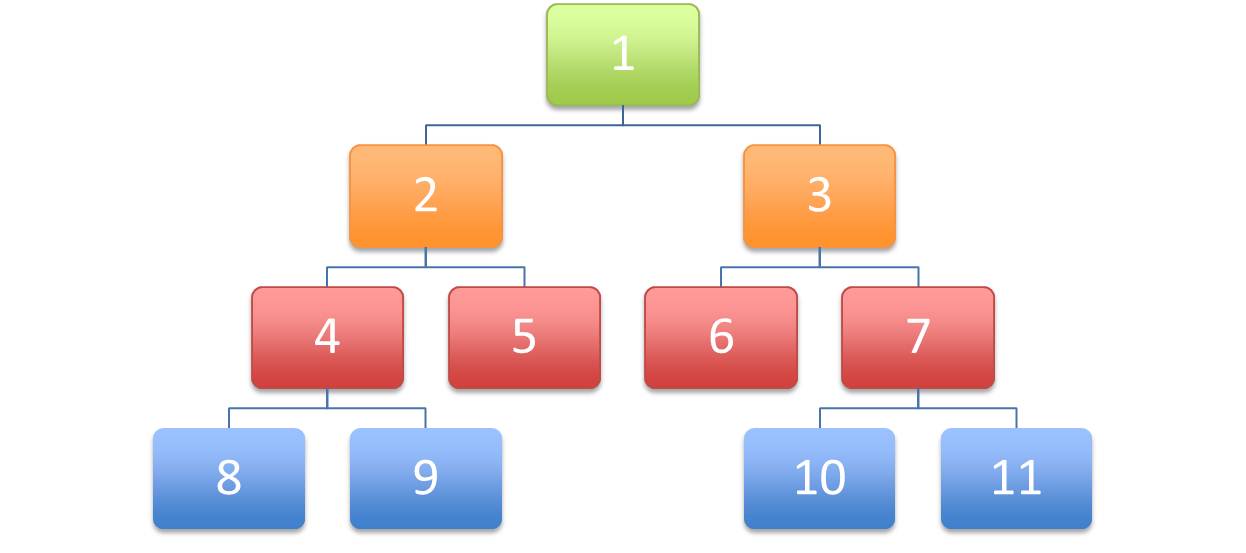}
\caption{Partition generated by METIS for system including dense layers
at two ends. \label{fig:partition_2}}
\end{figure}

In many cases, self-energy functions add two $N_{x}\times N_{x}$
dense blocks into the input sparse matrix $\mathbf{A}$ in first and
last block (see Figure \ref{fig:sparsity_A_H}). One possible partition
combines the two contacts together with the middle separator 1,
as shown in Figure \ref{fig:partition_1}.
The weakness of such clustering is the size of the first separator
(or root region). A larger separator increases the computational cost
spent on this level. More descendant blocks are coupled with this
separator and the total number of operations will increase dramatically.
Another partition that satisfies the previous two rules is plotted
in Figure \ref{fig:partition_2}.

\section{Numerical Experiments}

This section describes numerical experiments on two simple models:
a super-lattice structure and a graphene nanotube.
The algorithm is implemented in a C code that is interfaced with
METIS \cite{Karypis:1998:FHQ:305219.305248}.

\subsection{Cost analysis}

First the complexity of the HSC extension is compared numerically
to the complexity of RGF.
A model device is considered where the system Hamiltonian
is discretized with a five-point stencil.
The left and right contact self-energies are neglected for this section only.
A typical partition is plotted in Figure \ref{fig:multilevel_device}.

The numerical estimate tracks the operation counts step by step
for all the matrix multiplications and matrix inversions throughout the code.
For a multiplication of two matrices with dimensions $i \times j$ and $j \times k$,
a total of $ijk$ operations is added.
For inversion of a matrix of dimension $i \times i$, $i^3$ operations are counted.

\begin{figure}[htbp]
\centering \includegraphics[width=0.7\textwidth]{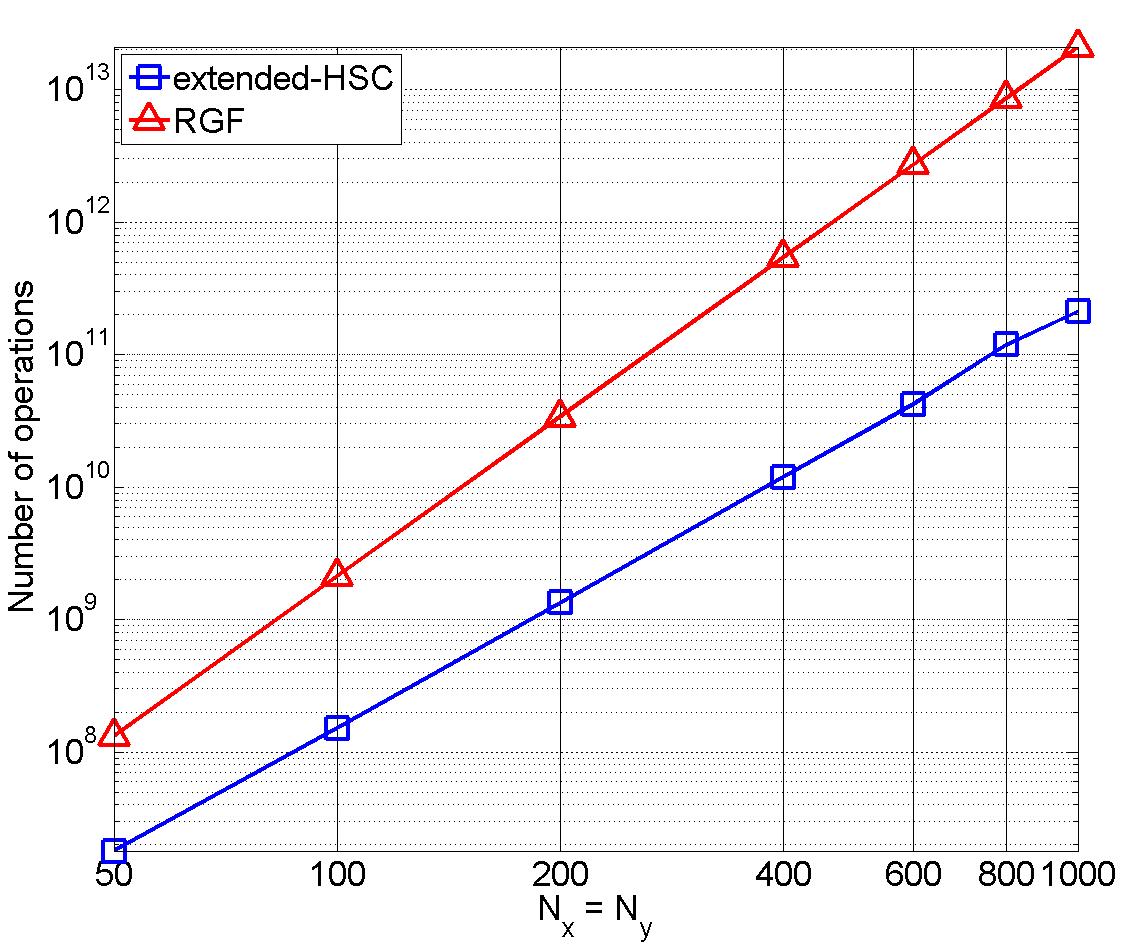}
\caption{Numerical count comparison for our algorithm (blue) and RGF (red).
\label{fig:cost}}
\end{figure}

Figure \ref{fig:cost} shows the cost comparison between the HSC extension
and RGF for two-dimensional square systems with the same number of
grid points (or atoms) per direction, {\it i.e. $N_{x}=N_{y}=N$}. 
The plot in logarithmic scale indicates that RGF exhibits 
a complexity of $\mathcal{O}(N^4)$, 
while the HSC-extension shows a $\mathcal{O}(N^3)$ complexity.

\subsection{Results}

\subsubsection{Super-lattice device}

\begin{figure}[htbp]
\centering \includegraphics[width=0.7\textwidth]{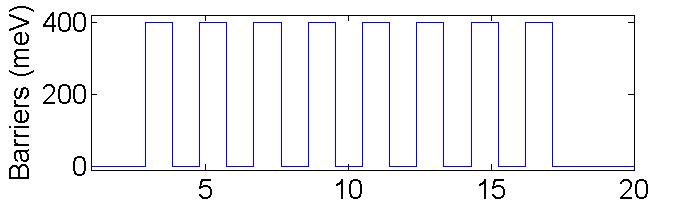}
\caption{Barrier structure for a model super-lattice device.}
\label{fig:modelSL}
\end{figure}

A super-lattice device is typically a multi-layered energy barriers system.
The device is composed of repeating junctions of energy barriers and wells. 
To verify the simulation results, a two-dimensional system of lengths 
$l_x = 25 \mathrm{ nm}$ and $l_y = 20 \mathrm{ nm}$ is considered and plotted in Figure \ref{fig:modelSL}.
Here the structure has eight barriers, each of width $1\mathrm{nm}$ 
and of height 400meV. 
The wells have a width of 1nm. 
The length of the left flat band region is $2\mathrm{nm}$ and the right flat band region is 3nm long.

\begin{figure}[htbp]
\centering \includegraphics[width=0.7\textwidth]{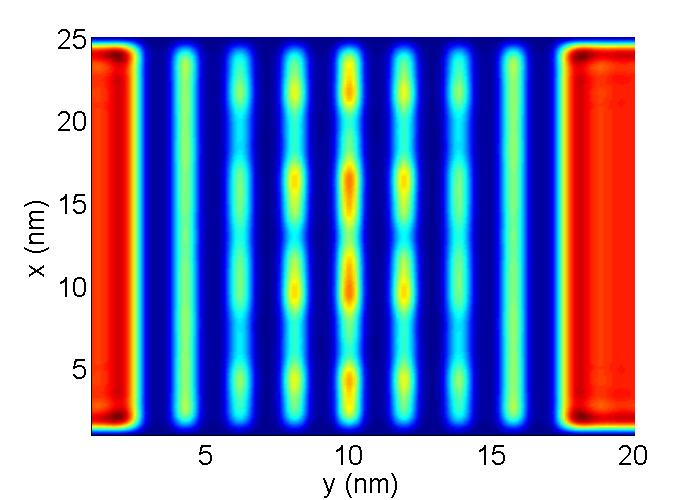}
\caption{Electron density for a model super-lattice device.}
\label{fig:modelSL_plot}
\end{figure}

\begin{figure}[htbp]
\centering \includegraphics[width=0.7\textwidth]{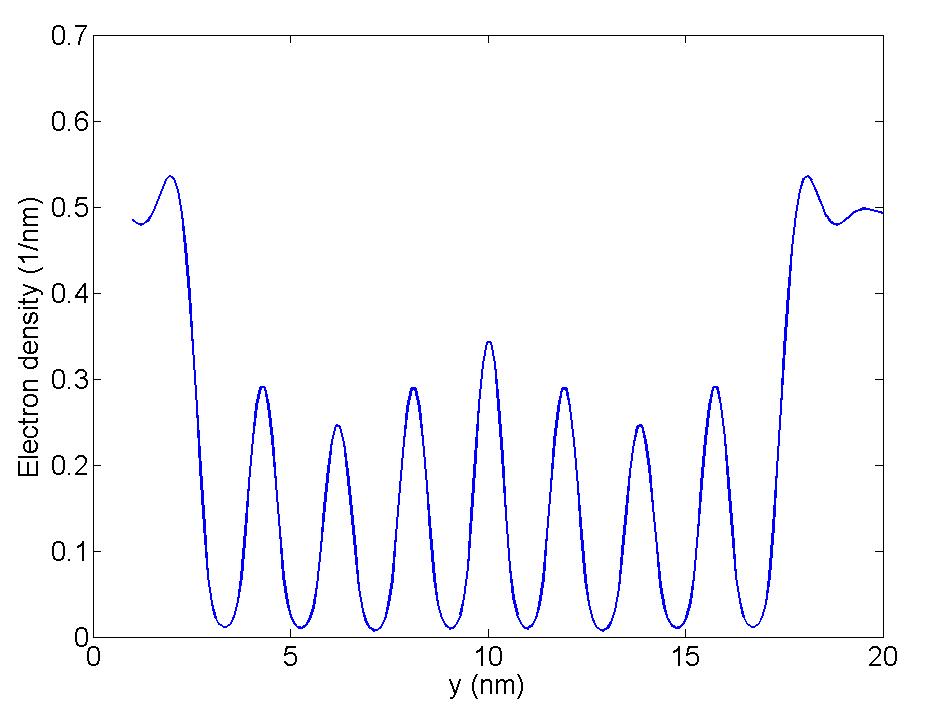}
\caption{Electron density in $y$ direction of the example super-lattice structure.}
\label{fig:rtddensy}
\end{figure}

A simulation with a five-point stencil discretization on a grid with 
spacing $\mathrm{d}x = \mathrm{d}y = 0.1 \mathrm{ nm}$
is made for the Fermi energy $E_f=140\mathrm{meV}$ and 500 energy points
uniformly distributed between 0 to $500\mathrm{meV}$.
The density of states, electron density, and current are calculated by the RGF method and
the extended-HSC method.
The output electron density is plotted in Figure \ref{fig:modelSL_plot}.
Linear electron density in the $y$-direction is illustrated in Figure \ref{fig:rtddensy}.
The figures indicate that the charge distribution in the barrier-well multi-layer junctions are
symmetric, as expected.

Next devices of lengths $l_x = N_x \times 0.1 \mathrm{nm}$ and
$l_y = N_y \times 0.1 \mathrm{nm}$ are used to compare the two algorithms.
The number of barriers is kept at 8.
The lengths for the two-sides flat region are adjusted according to
the lengths of device $l_x$ and $l_y$.
The other parameters remain unchanged.

\begin{figure}[htbp]
\centering \includegraphics[width=0.9\textwidth]{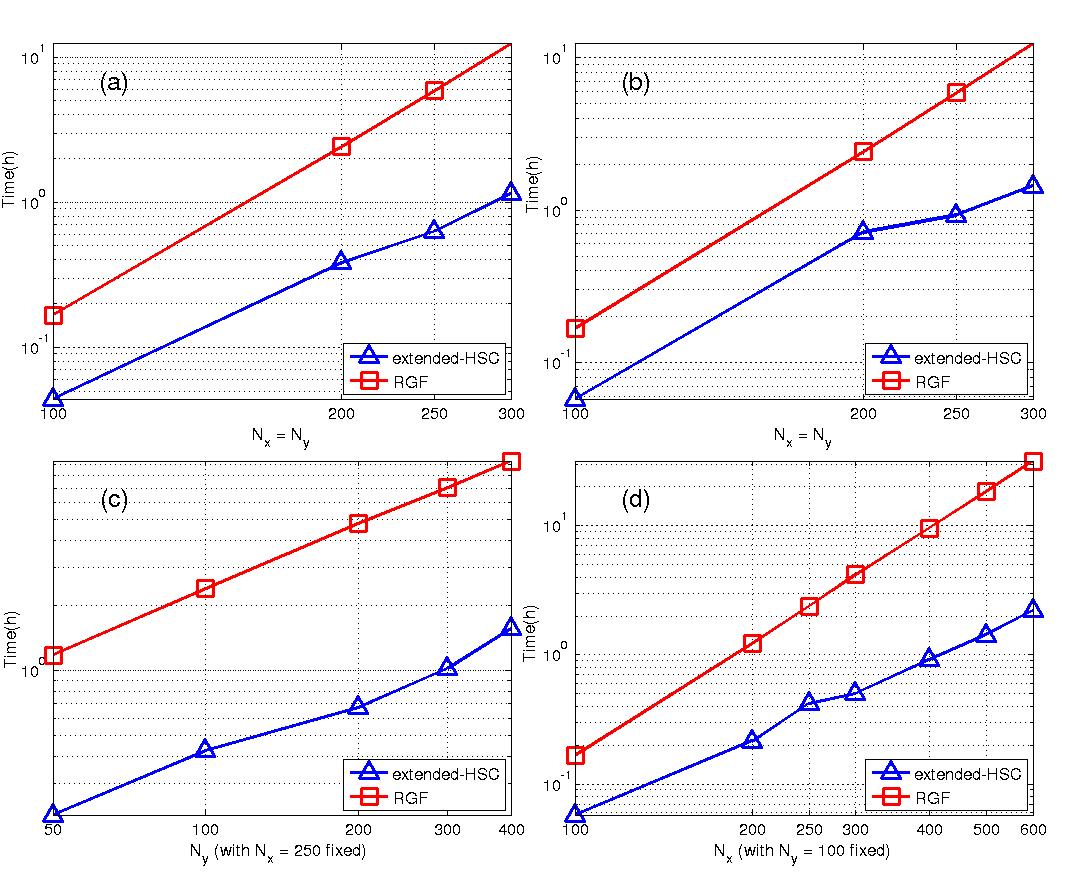}
\caption{Superlattice device NEGF simulation computation time comparison for
RGF and our methods, all systems grid spacing is 0.1nm. (a) Square
system of with diagonal self-energy matrix; (b) Square system of
with dense self-energy matrix; (c) For systems in this plot, the
length in the $x$-direction is fixed at 25nm while the length
in the $y$-direction is increased.
(d) For systems in this plot, the length in the $y$-direction is fixed at
10nm while the length in the $x$-direction is increased.
Dense self-energy matrices are used in (c) and (d) devices.}
\label{fig:resultSL}
\end{figure}

In Figure \ref{fig:resultSL}(a), diagonal self-energy matrices are used
for the left and right contacts.
Calculation times are compared for square systems--- 
{\it i.e. $N_x = N_y=N$} ---
and plotted in Figure \ref{fig:resultSL}(a).
As expected, the HSC-extension exhibits smaller CPU times and a 
complexity of $\mathcal{O}(N^3)$ while RGF's complexity is 
$\mathcal{O}(N^4)$.

In Figure \ref{fig:resultSL}(b), CPU times for square systems with
dense self-energy matrices for both contacts are plotted.
Here again the HSC-extension exhibits smaller CPU times.
A complexity $\mathcal{O}(N^3)$ for HSC-extension compared with 
$\mathcal{O}(N^4)$ for RGF can be seen.

Figure \ref{fig:resultSL}(c) plots CPU times for rectangular devices where
$l_x = 25 \mathrm{nm}$ (or $N_x = 250$) and the length in the $y$-direction
is varied.
Dense self-energy blocks for the left and right contacts are employed.
The implementation of RGF is biased towards the $x$-direction so that
its complexity is $\mathcal{O}\left(N_y \right)$.
The linear trend is clearly present in the plot.
For the HSC-extension, similarly to the cost of a sparse $\mathbf{LDL}^T$ factorization,
the computational cost is $\mathcal{O}\left(N_y \right)$.
Clearly, the constant for RGF is larger for this device.

To illustrate the dependence of this constant with respect to $N_x$,
Figure \ref{fig:resultSL}(d) plots the CPU times when $N_y$ is fixed
and $N_x$ is varied.
The recorded CPU times illustrate that the RGF method has
an asymptotic complexity $\mathcal{O}\left( N_x^3 \right)$,
while the HSC extension exhibits a complexity $\mathcal{O}\left( N_x^2 \right)$.
So, for rectangular devices, the RGF method has a complexity
$\mathcal{O}(N_x^3 N_y)$ and the HSC-extension a complexity 
$\mathcal{O}(N_x^2 N_y)$.

\subsubsection{Graphene}

Graphene is one of the most promising next-generation materials. Its remarkable electric properties, such as high carrier mobility and zero band gaps, generate a rapidly increasing interest in the electronic device community. Since 2007, many advances in graphene-based transistor development have been reported. \cite{Schwierz:2010fk}

\begin{figure}[htbp]
\begin{center}
\includegraphics[width=9cm]{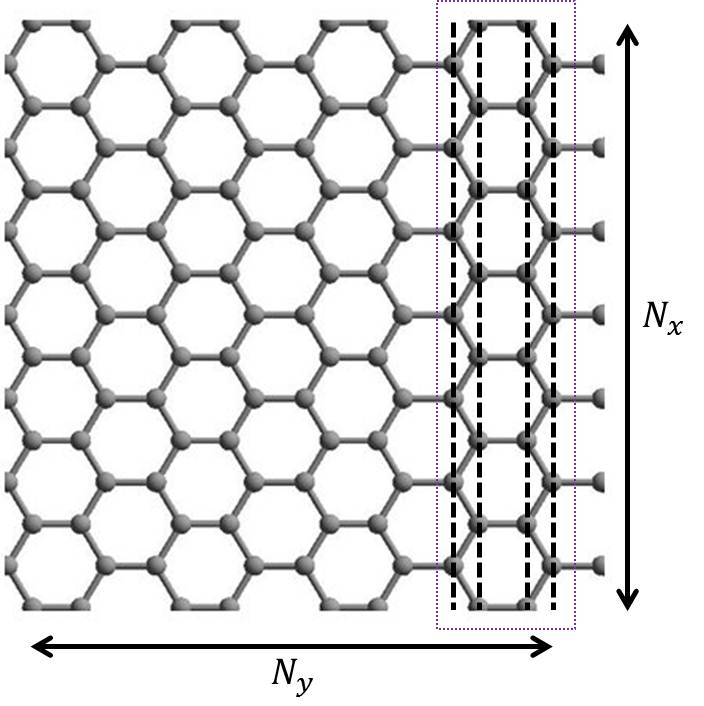}
\end{center}
\caption{Graphene hexagonal structure decomposed by tight biding method.
Dashed rectangular illustrates one repeating hexagon layer.
Dashed lines represent inner four atom layers, showing the atoms ordering in tight binding
Hamiltonian construction.}
\label{fig:grap}
\end{figure}

The NEGF simulation of graphene transport is based on tight binding method, which yields a four-point-stencil Hamiltonian due to system decomposition of carbon atoms coupling (see Figure \ref{fig:grap}). In the numerical experiments, armchair planar graphene nanoribbon structures are simulated. The on-site energy for each carbon atom is 0 and the hopping parameter between two nearest carbon atoms is -3.1eV. The Fermi energy is set to 0. The simulation is run for only one energy point $E=0.5$eV.

\begin{figure}[htbp]
\begin{center}
\includegraphics[width=10cm]{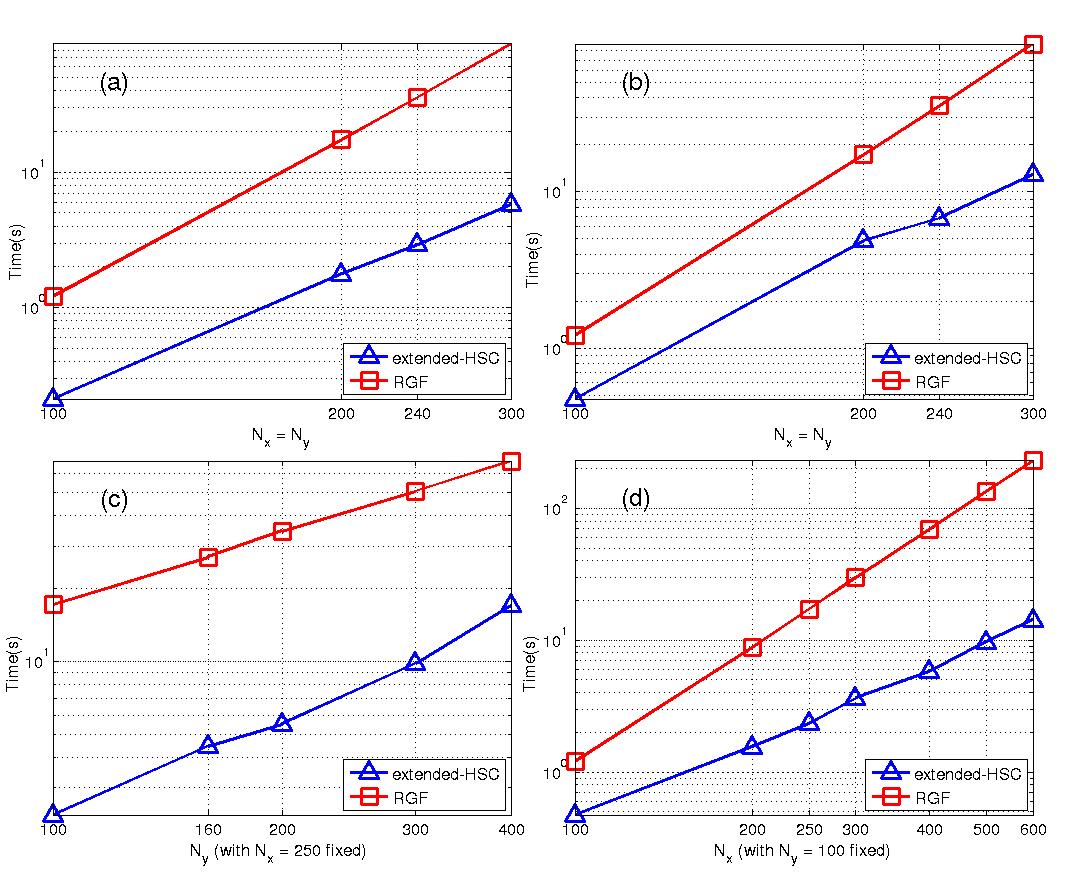}
\end{center}
\caption{Graphene device NEGF simulation computation time comparison for
RGF and the HSC extension, based on tight binding theory.
(a) Square system of with diagonal self-energy matrix.
(b) Square system of with dense self-energy matrix.
(c) For systems in this plot, the number of atoms in the $x$-direction is set to
$N_x = 250$, while the number of atoms in the $y$-direction
is increased.
(d) For systems in this plot, the number of atoms in the $y$-direction is set to
$N_y = 100$, while the number of atoms in the $x$-direction
is increased.
Dense self-energy matrices are used in (c) and (d) devices.}
\label{fig:resultGP}
\end{figure}

Simulation timings are plotted in Figure \ref{fig:resultGP} for graphene structures of different sizes. To minimize the dimension of blocks to invert in RGF, one layer of hexagonal structure is divided into four layers (see the dashed lines in Figure \ref{fig:grap}). Conclusions on the asymptotic complexity remain unchanged. Namely, the HSC extension is more efficient than the RGF method for square and rectangular structures. The complexity of RGF for four-point stencil behaves as $\mathcal{O}(N_x^3N_y)$ with the same constant as five-point stencil, which is expected due to the layered system partition. In Figure \ref{fig:resultGP}(a) and (b), for a four-point stencil system, our HSC-extension exhibits a complexity of $\mathcal{O}(N^3)$ for square system. 
Figure \ref{fig:resultGP}(c) and (d) also demonstrate a complexity growing linearly with $N_y$ and, respectively quadratically with $N_x$, when $N_x$, respectively $N_y$, is fixed. The final result is a complexity $\mathcal{O}(N_x^2N_y)$. The constant in front of $N_x^2N_y$ for the extended HSC method is smaller for the graphene structures than for the superlattice devices. This reduction is explained by a more efficient partitioning of four-point stencil systems by METIS .

\section{Conclusion}

In this paper, an approach to calculate the charge density
in nanoscale devices, within the context of the non equilibrium Green's
function approach, is presented.
This work exploits recent advances to use an established graph partitioning
method, namely the nested dissection method.
This contribution does not require any processing of the partition and
it can handle open boundary conditions, represented by full self-energy
matrices.
The key ingredients are an efficient sparse block $\mathbf{LDL}^T$-factorization
and an appropriate order of operations to preserve the sparsity as much as
possible.
The resulting algorithm was illustrated on a quantum well superlattice
and a carbon nanotube, which are represented by a continuum and tight
binding Hamiltonian respectively, and demonstrated a significant speed up
over the recursive method RGF.

\section*{Acknowledgements}

The authors acknowledge the support by the National Science Foundation 
under Grant ECCS-1231927.
The authors thank the anonymous referees
for comments that led to improvements of the manuscript.

\bibliographystyle{plain}
\bibliography{Refs_GrGless}

\appendix

\section{Description of the Algorithm for a Three-level Tree}

\begin{figure}[htbp]
\centering \includegraphics[width=0.45\textwidth]{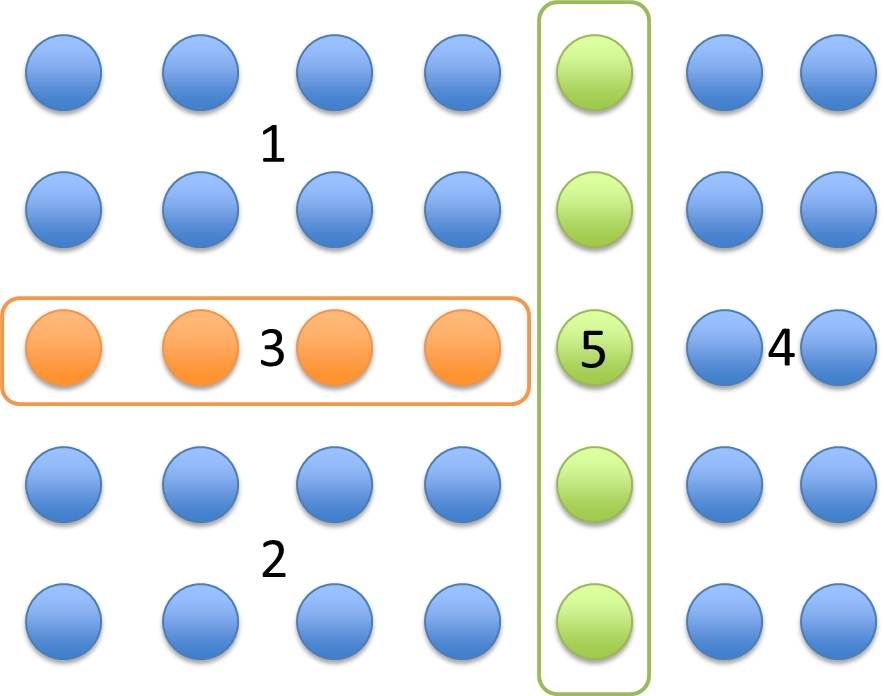}  \includegraphics[width=0.45\textwidth]{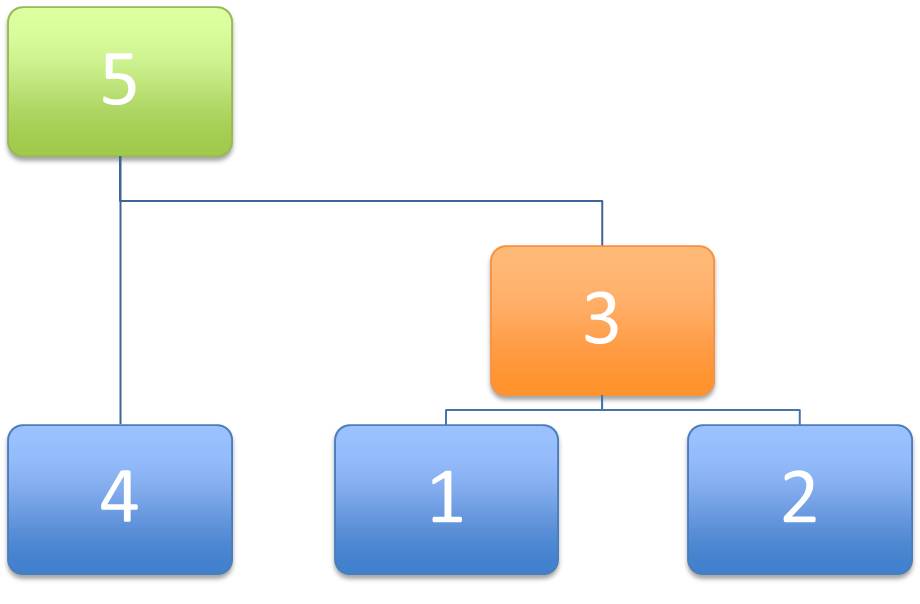}
\caption{Partition for a three-level system.\label{fig:partition_3level}}
\end{figure}

In order to make the extension more comprehensive, a description of
the HSC extension is given for a three-level system (see Figure \ref{fig:partition_3level}).

 The size of the partition is chosen to make the description as relevant
as possible without becoming overcomplicated. The first level contains
regions 1, 2, and 4. The second-level refers to region 3 and the top
level is the root or region 5. To illustrate the algorithm, steps
from Section \ref{sub:description_multi} are described for this particular
device.

When a five-point stencil is used for discretization, the structure
of the matrix $\mathbf{A}$, after re-ordering, is
\[
\mathbf{A}=\left[\begin{array}{ccccc}
\mathbf{A}_{11} & \mathbf{0} & \mathbf{A}_{13} & \mathbf{0} & \mathbf{A}_{15}\\
\mathbf{0} & \mathbf{A}_{22} & \mathbf{A}_{23} & \mathbf{0} & \mathbf{A}_{25}\\
\mathbf{A}_{13}^{T} & \mathbf{A}_{23}^{T} & \mathbf{A}_{33} & \mathbf{0} & \mathbf{A}_{35}\\
\mathbf{0} & \mathbf{0} & \mathbf{0} & \mathbf{A}_{44} & \mathbf{A}_{45}\\
\mathbf{A}_{15}^{T} & \mathbf{A}_{25}^{T} & \mathbf{A}_{35}^{T} & \mathbf{A}_{45}^{T} & \mathbf{A}_{55}
\end{array}\right]
\]
The blocks $\mathbf{A}_{11}$ and $\mathbf{A}_{44}$ are dense for
representing the contacts with the semi-infinite leads.

\begin{figure}[htbp]
\centering \includegraphics[width=0.5\textwidth]{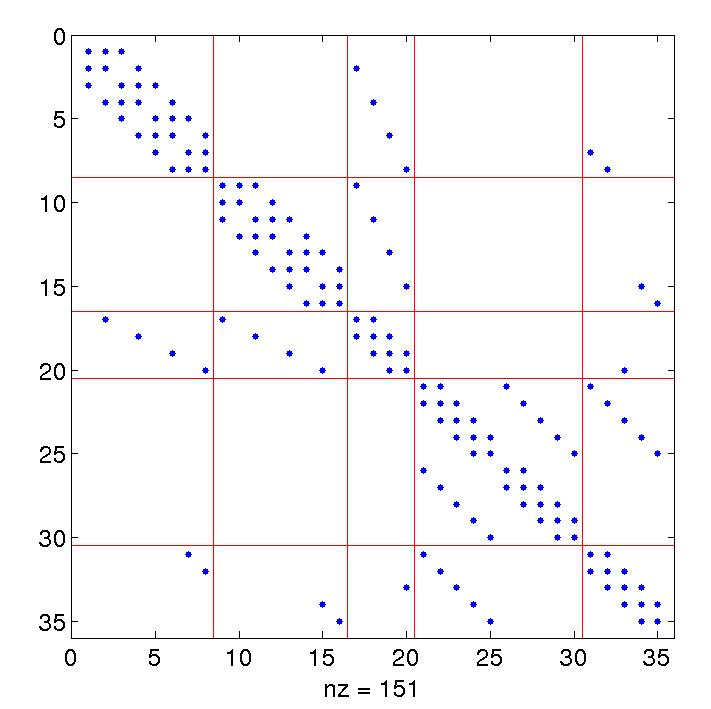}\includegraphics[width=0.5\textwidth]{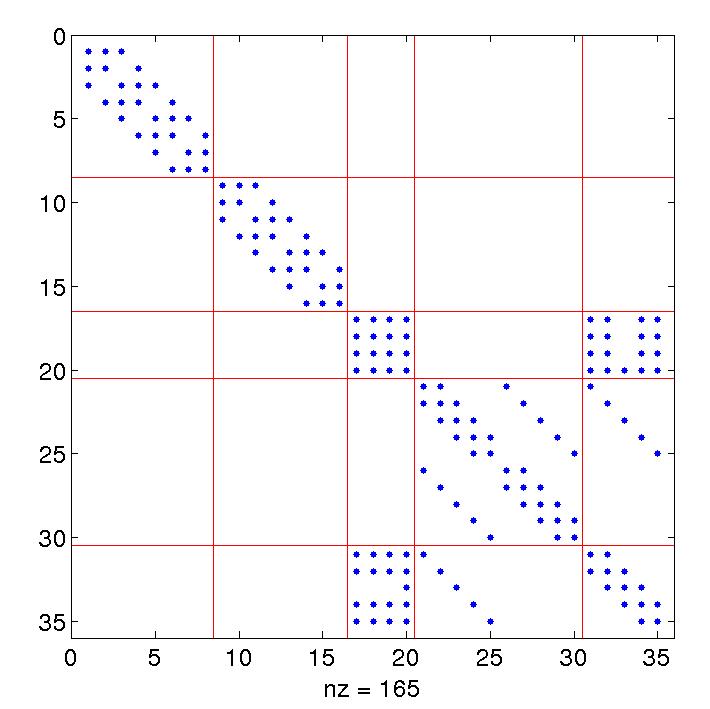}
\caption{Sparsity of matrix $\mathbf{A}$ and matrix $\mathbf{A}^{(1)}$.\label{fig:sparsity_A0_A1}}
\end{figure}

Recall that $\mathbf{A}^{(0)}$ is equal to $\mathbf{A}$. Then inner
points in regions 1, 2, and 4 are eliminated by block Gaussian elimination
--- the effects of the inner points in regions 1, 2, and 4 are folded
over their boundary. This first step yields the matrix $\mathbf{A}^{(1)}$
\[
\mathbf{A}^{(1)}=\left[\begin{array}{ccccc}
\mathbf{A}_{11} & \mathbf{0} & \mathbf{0} & \mathbf{0} & \mathbf{0}\\
\mathbf{0} & \mathbf{A}_{22} & \mathbf{0} & \mathbf{0} & \mathbf{0}\\
\mathbf{0} & \mathbf{0} & \mathbf{A}_{33}^{(1)} & \mathbf{0} & \mathbf{A}_{35}^{(1)}\\
\mathbf{0} & \mathbf{0} & \mathbf{0} & \mathbf{A}_{44} & \mathbf{0}\\
\mathbf{0} & \mathbf{0} & \left(\mathbf{A}_{35}^{(1)}\right)^{T} & \mathbf{0} & \mathbf{A}_{55}^{(1)}
\end{array}\right]
\]
where the updated matrices are
\begin{alignat*}{1}
\mathbf{A}_{33}^{(1)} & =\mathbf{A}_{33}-\mathbf{A}_{13}^{T}\mathbf{A}_{11}^{-1}\mathbf{A}_{13}-\mathbf{A}_{23}^{T}\mathbf{A}_{22}^{-1}\mathbf{A}_{23}
,
\\
\mathbf{A}_{35}^{(1)} & =\mathbf{A}_{35}-\mathbf{A}_{13}^{T}\mathbf{A}_{11}^{-1}\mathbf{A}_{15}-\mathbf{A}_{23}^{T}\mathbf{A}_{22}^{-1}\mathbf{A}_{25}
,
\\
\mathbf{A}_{55}^{(1)} & =\mathbf{A}_{55}-\mathbf{A}_{15}^{T}\mathbf{A}_{11}^{-1}\mathbf{A}_{15}-\mathbf{A}_{25}^{T}\mathbf{A}_{22}^{-1}\mathbf{A}_{25}-\mathbf{A}_{45}^{T}\mathbf{A}_{44}^{-1}\mathbf{A}_{45}
.
\end{alignat*}
After folding the effects of regions 1 and 2, block $\mathbf{A}_{35}^{(1)}$
is now a dense block. Figure \ref{fig:sparsity_A0_A1} illustrates
the change of sparsity between $\mathbf{A}$ and $\mathbf{A}^{(1)}$.

In the next step, the remaining off-diagonal blocks are eliminated
to obtain the matrix $\mathbf{A}^{(2)}$,
\[
\mathbf{A}^{(2)}=\left[\begin{array}{ccccc}
\mathbf{A}_{11} & \mathbf{0} & \mathbf{0} & \mathbf{0} & \mathbf{0}\\
\mathbf{0} & \mathbf{A}_{22} & \mathbf{0} & \mathbf{0} & \mathbf{0}\\
\mathbf{0} & \mathbf{0} & \mathbf{A}_{33}^{(1)} & \mathbf{0} & \mathbf{0}\\
\mathbf{0} & \mathbf{0} & \mathbf{0} & \mathbf{A}_{44} & \mathbf{0}\\
\mathbf{0} & \mathbf{0} & \mathbf{0} & \mathbf{0} & \mathbf{A}_{55}^{(2)}
\end{array}\right]
\]
where the block $\mathbf{A}_{55}^{(2)}$ is
\[
\mathbf{A}_{55}^{(2)}=\mathbf{A}_{55}^{(1)}-\left(\mathbf{A}_{35}^{(1)}\right)^{T}\left(\mathbf{A}_{33}^{(1)}\right)^{-1}\mathbf{A}_{35}^{(1)}.
\]
Note that several blocks are unchanged, like $\mathbf{A}_{11}$, 
$\mathbf{A}_{22}$, $\mathbf{A}_{33}^{(1)}$, and $\mathbf{A}_{44}$.
 The next step is written as the inversion of $\mathbf{A}^{(2)}$,
which is symmetric and block diagonal,
\[
\mathbf{G}^{(2)}=\left[\begin{array}{ccccc}
\mathbf{A}_{11}^{-1} & \mathbf{0} & \mathbf{0} & \mathbf{0} & \mathbf{0}\\
\mathbf{0} & \mathbf{A}_{22}^{-1} & \mathbf{0} & \mathbf{0} & \mathbf{0}\\
\mathbf{0} & \mathbf{0} & \left(\mathbf{A}_{33}^{(1)}\right)^{-1} & \mathbf{0} & \mathbf{0}\\
\mathbf{0} & \mathbf{0} & \mathbf{0} & \mathbf{A}_{44}^{-1} & \mathbf{0}\\
\mathbf{0} & \mathbf{0} & \mathbf{0} & \mathbf{0} & \left(\mathbf{A}_{55}^{(2)}\right)^{-1}
\end{array}\right].
\]
This operation requires only the inversion of the block $\mathbf{A}_{55}^{(2)}$.
All the other blocks have been inverted during the folding steps.

\begin{figure}[htbp]
\centering \includegraphics[width=0.5\textwidth]{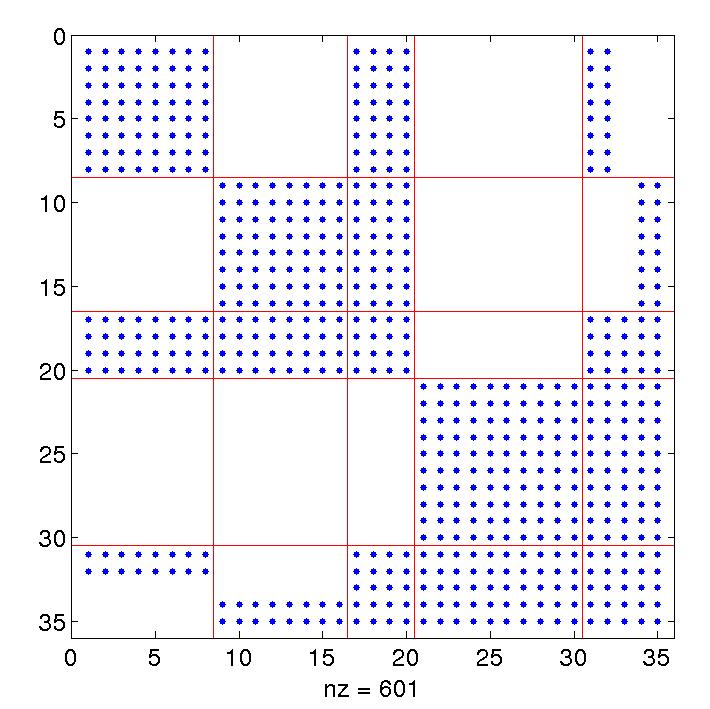}
\caption{Sparsity of matrix $\mathbf{G}^{(0)}$.\label{fig:sparsity_G0}}
\end{figure}

Next diagonal blocks of $\mathbf{G}^{r}$ are extracted one level
at a time. Starting from the main root (or separator), blocks at level
2 are updated to obtain
\[
\mathbf{G}^{(1)}=\left[\begin{array}{ccccc}
\mathbf{G}_{11}^{(2)} & \mathbf{0} & \mathbf{0} & \mathbf{0} & \mathbf{0}\\
\mathbf{0} & \mathbf{G}_{22}^{(2)} & \mathbf{0} & \mathbf{0} & \mathbf{0}\\
\mathbf{0} & \mathbf{0} & \mathbf{G}_{33}^{(1)} & \mathbf{0} & \mathbf{G}_{35}^{(1)}\\
\mathbf{0} & \mathbf{0} & \mathbf{0} & \mathbf{G}_{44}^{(2)} & \mathbf{0}\\
\mathbf{0} & \mathbf{0} & \left(\mathbf{G}_{35}^{(1)}\right)^{T} & \mathbf{0} & \mathbf{G}_{55}^{(2)}
\end{array}\right]
\]
with
\begin{alignat*}{1}
\mathbf{G}_{35}^{(1)} & =-\left(\mathbf{A}_{33}^{(1)}\right)^{-1}\mathbf{A}_{35}^{(1)}\mathbf{G}_{55}^{(2)}=\boldsymbol{\Psi}_{35}\mathbf{G}_{55}^{(2)}
,
\\
\mathbf{G}_{33}^{(1)} & =\mathbf{G}_{33}^{(2)}-\left(\mathbf{A}_{33}^{(1)}\right)^{-1}\mathbf{A}_{35}^{(1)}\left(\mathbf{G}_{35}^{(1)}\right)^{T}=\mathbf{G}_{33}^{(2)}+\boldsymbol{\Psi}_{35}\left(\mathbf{G}_{35}^{(1)}\right)^{T}
.
\end{alignat*}
Finally, blocks for regions 1, 2, and 4 are updated, yielding the
matrix $\mathbf{G}^{(0)}$,
\[
\mathbf{G}^{(0)}=\left[\begin{array}{ccccc}
\mathbf{G}_{11}^{(0)} & \mathbf{0} & \mathbf{G}_{13}^{(0)} & \mathbf{0} & \mathbf{G}_{15}^{(0)}\\
\mathbf{0} & \mathbf{G}_{22}^{(0)} & \mathbf{G}_{23}^{(0)} & \mathbf{0} & \mathbf{G}_{25}^{(0)}\\
\left(\mathbf{G}_{13}^{(0)}\right)^{T} & \left(\mathbf{G}_{23}^{(0)}\right)^{T} & \mathbf{G}_{33}^{(1)} & \mathbf{0} & \mathbf{G}_{35}^{(1)}\\
\mathbf{0} & \mathbf{0} & \mathbf{0} & \mathbf{G}_{44}^{(0)} & \mathbf{G}_{45}^{(0)}\\
\left(\mathbf{G}_{15}^{(0)}\right)^{T} & \left(\mathbf{G}_{25}^{(0)}\right)^{T} & \left(\mathbf{G}_{35}^{(1)}\right)^{T} & \left(\mathbf{G}_{45}^{(0)}\right)^{T} & \mathbf{G}_{55}^{(2)}
\end{array}\right].
\]
Blocks for region 4 are satisfying
\begin{alignat*}{1}
\mathbf{G}_{45}^{(0)} & =-\left(\mathbf{A}_{44}^{(0)}\right)^{-1}\mathbf{A}_{45}^{(0)}\mathbf{G}_{55}^{(2)}=\boldsymbol{\Psi}_{45}\mathbf{G}_{55}^{(2)}\\
\mathbf{G}_{44}^{(0)} & =\mathbf{G}_{44}^{(2)}-\left(\mathbf{A}_{44}^{(0)}\right)^{-1}\mathbf{A}_{45}^{(0)}\left(\mathbf{G}_{45}^{(0)}\right)^{T}
\end{alignat*}
Blocks for region 1 are defined by
\begin{alignat*}{1}
\mathbf{G}_{15}^{(0)} & =-\left(\mathbf{A}_{11}^{(0)}\right)^{-1}\mathbf{A}_{15}^{(0)}\mathbf{G}_{55}^{(2)}-\left(\mathbf{A}_{11}^{(0)}\right)^{-1}\mathbf{A}_{13}^{(0)}\mathbf{G}_{35}^{(1)}\\
\mathbf{G}_{13}^{(0)} & =-\left(\mathbf{A}_{11}^{(0)}\right)^{-1}\mathbf{A}_{13}^{(0)}\mathbf{G}_{33}^{(1)}-\left(\mathbf{A}_{11}^{(0)}\right)^{-1}\mathbf{A}_{15}^{(0)}\left(\mathbf{G}_{35}^{(1)}\right)^{T}\\
\mathbf{G}_{11}^{(0)} & =\left(\mathbf{A}_{11}^{(0)}\right)^{-1}-\left(\mathbf{A}_{11}^{(0)}\right)^{-1}\mathbf{A}_{13}^{(0)}\left(\mathbf{G}_{13}^{(0)}\right)^{T}-\left(\mathbf{A}_{11}^{(0)}\right)^{-1}\mathbf{A}_{15}^{(0)}\left(\mathbf{G}_{15}^{(0)}\right)^{T}
\end{alignat*}
and blocks for region 2
\begin{alignat*}{1}
\mathbf{G}_{25}^{(0)} & =-\left(\mathbf{A}_{22}^{(0)}\right)^{-1}\mathbf{A}_{25}^{(0)}\mathbf{G}_{55}^{(2)}-\left(\mathbf{A}_{22}^{(0)}\right)^{-1}\mathbf{A}_{23}^{(0)}\mathbf{G}_{35}^{(1)}\\
\mathbf{G}_{23}^{(0)} & =-\left(\mathbf{A}_{22}^{(0)}\right)^{-1}\mathbf{A}_{23}^{(0)}\mathbf{G}_{33}^{(1)}-\left(\mathbf{A}_{22}^{(0)}\right)^{-1}\mathbf{A}_{25}^{(0)}\left(\mathbf{G}_{35}^{(1)}\right)^{T}\\
\mathbf{G}_{22}^{(0)} & =\left(\mathbf{A}_{22}^{(0)}\right)^{-1}-\left(\mathbf{A}_{22}^{(0)}\right)^{-1}\mathbf{A}_{23}^{(0)}\left(\mathbf{G}_{23}^{(0)}\right)^{T}-\left(\mathbf{A}_{22}^{(0)}\right)^{-1}\mathbf{A}_{25}^{(0)}\left(\mathbf{G}_{25}^{(0)}\right)^{T}
\end{alignat*}
Figure \ref{fig:sparsity_G0} displays the sparsity of the resulting
matrix $\mathbf{G}^{(0)}$.
 All the entries in $\mathbf{G}^{(0)}$ are equal to their corresponding
entries in $\mathbf{G}^{r}$.

The computation of diagonal blocks in $\mathbf{G}^{<}$ are described
for the same device. First the matrix 
$\mathbf{N}$ is computed,
\[
\mathbf{N}
=
\boldsymbol{\Sigma}^{<} \left( \mathbf{G}^{<} \right)^\dagger
=
\left[\begin{array}{ccccc}
\boldsymbol{\Sigma}_{11}^{<}\overline{\mathbf{G}_{11}^{(0)}} & \mathbf{0} & \boldsymbol{\Sigma}_{11}^{<}\overline{\mathbf{G}_{13}^{(0)}} & \mathbf{0} & \boldsymbol{\Sigma}_{11}^{<}\overline{\mathbf{G}_{15}^{(0)}}\\
\mathbf{0} & \boldsymbol{\Sigma}_{22}^{<}\overline{\mathbf{G}_{22}^{(0)}} & \boldsymbol{\Sigma}_{22}^{<}\overline{\mathbf{G}_{23}^{(0)}} & \mathbf{0} & \boldsymbol{\Sigma}_{22}^{<}\overline{\mathbf{G}_{25}^{(0)}}\\
\boldsymbol{\Sigma}_{33}^{<}\left(\mathbf{G}_{13}^{(0)}\right)^{\dagger} & \boldsymbol{\Sigma}_{33}^{<}\left(\mathbf{G}_{23}^{(0)}\right)^{\dagger} & \boldsymbol{\Sigma}_{33}^{<}\overline{\mathbf{G}_{33}^{(0)}} & \mathbf{0} & \boldsymbol{\Sigma}_{33}^{<}\overline{\mathbf{G}_{35}^{(0)}}\\
\mathbf{0} & \mathbf{0} & \mathbf{0} & \boldsymbol{\Sigma}_{44}^{<}\overline{\mathbf{G}_{44}^{(0)}} & \boldsymbol{\Sigma}_{44}^{<}\overline{\mathbf{G}_{45}^{(0)}}\\
\boldsymbol{\Sigma}_{55}^{<}\left(\mathbf{G}_{15}^{(0)}\right)^{\dagger} & \boldsymbol{\Sigma}_{55}^{<}\left(\mathbf{G}_{25}^{(0)}\right)^{\dagger} & \boldsymbol{\Sigma}_{55}^{<}\left(\mathbf{G}_{35}^{(0)}\right)^{\dagger} & \boldsymbol{\Sigma}_{55}^{<}\left(\mathbf{G}_{45}^{(0)}\right)^{\dagger} & \boldsymbol{\Sigma}_{55}^{<}\overline{\mathbf{G}_{55}^{(0)}}
\end{array}\right].
\]
Set $\mathbf{N}^{(0)}=\mathbf{N}$. Next the lower level clusters
are folded into the higher ones to obtain the matrix $\mathbf{N}^{(1)}$
\[
\mathbf{N}^{(1)}=\left[\begin{array}{ccccc}
\mathbf{N}_{11}^{(0)} & \mathbf{0} & \mathbf{N}_{13}^{(0)} & \mathbf{0} & \mathbf{N}_{15}^{(0)}\\
\mathbf{0} & \mathbf{N}_{22}^{(0)} & \mathbf{N}_{23}^{(0)} & \mathbf{0} & \mathbf{N}_{25}^{(0)}\\
\mathbf{N}_{31}^{(0)} & \mathbf{N}_{32}^{(0)} & \mathbf{N}_{33}^{(1)} & \mathbf{0} & \mathbf{N}_{35}^{(1)}\\
\mathbf{0} & \mathbf{0} & \mathbf{0} & \mathbf{N}_{44}^{(0)} & \mathbf{N}_{45}^{(0)}\\
\mathbf{N}_{51}^{(0)} & \mathbf{N}_{52}^{(0)} & \mathbf{N}_{53}^{(1)} & \mathbf{N}_{54}^{(0)} & \mathbf{N}_{55}^{(1)}
\end{array}\right]
\]
where the updated blocks are
\begin{alignat*}{1}
\mathbf{N}_{33}^{(1)} & =\mathbf{N}_{33}^{(0)}-\left(\mathbf{A}_{13}\right)^{T}\mathbf{A}_{11}^{-1}\mathbf{N}_{13}^{(0)}-\left(\mathbf{A}_{23}\right)^{T}\mathbf{A}_{22}^{-1}\mathbf{N}_{23}^{(0)}\\
\mathbf{N}_{55}^{(1)} & =\mathbf{N}_{55}^{(0)}-\left(\mathbf{A}_{15}\right)^{T}\mathbf{A}_{11}^{-1}\mathbf{N}_{15}^{(0)}-\left(\mathbf{A}_{25}\right)^{T}\mathbf{A}_{22}^{-1}\mathbf{N}_{25}^{(0)}-\left(\mathbf{A}_{45}\right)^{T}\mathbf{A}_{44}^{-1}\mathbf{N}_{45}^{(0)}\\
\mathbf{N}_{35}^{(1)} & =\mathbf{N}_{35}^{(0)}-\left(\mathbf{A}_{13}\right)^{T}\mathbf{A}_{11}^{-1}\mathbf{N}_{15}^{(0)}-\left(\mathbf{A}_{23}\right)^{T}\mathbf{A}_{22}^{-1}\mathbf{N}_{25}^{(0)}\\
\mathbf{N}_{53}^{(1)} & =\mathbf{N}_{53}^{(0)}-\left(\mathbf{A}_{15}\right)^{T}\mathbf{A}_{11}^{-1}\mathbf{N}_{13}^{(0)}-\left(\mathbf{A}_{25}\right)^{T}\mathbf{A}_{22}^{-1}\mathbf{N}_{23}^{(0)}
\end{alignat*}
For the top level, the block for region 5 is updated
\[
\mathbf{N}^{(2)}=\left[\begin{array}{ccccc}
\mathbf{N}_{11}^{(0)} & \mathbf{0} & \mathbf{N}_{13}^{(0)} & \mathbf{0} & \mathbf{N}_{15}^{(0)}\\
\mathbf{0} & \mathbf{N}_{22}^{(0)} & \mathbf{N}_{23}^{(0)} & \mathbf{0} & \mathbf{N}_{25}^{(0)}\\
\mathbf{N}_{31}^{(0)} & \mathbf{N}_{32}^{(0)} & \mathbf{N}_{33}^{(1)} & \mathbf{0} & \mathbf{N}_{35}^{(1)}\\
\mathbf{0} & \mathbf{0} & \mathbf{0} & \mathbf{N}_{44}^{(0)} & \mathbf{N}_{45}^{(0)}\\
\mathbf{N}_{51}^{(0)} & \mathbf{N}_{52}^{(0)} & \mathbf{N}_{53}^{(1)} & \mathbf{N}_{54}^{(0)} & \mathbf{N}_{55}^{(2)}
\end{array}\right]
\]
with
\[
\mathbf{N}_{55}^{(2)}=\mathbf{N}_{55}^{(1)}-\left(\mathbf{A}_{35}^{(1)}\right)^{T}\left(\mathbf{A}_{33}^{(1)}\right)^{-1}\mathbf{N}_{35}^{(1)}.
\]

The next step is a block-diagonal multiplication
\[
\mathbf{P}^{(2)}=\left[\begin{array}{ccccc}
\mathbf{A}_{11}^{-1}\mathbf{N}_{11}^{(0)} & \mathbf{0} & \mathbf{A}_{11}^{-1}\mathbf{N}_{13}^{(0)} & \mathbf{0} & \mathbf{A}_{11}^{-1}\mathbf{N}_{15}^{(0)}\\
\mathbf{0} & \mathbf{A}_{22}^{-1}\mathbf{N}_{22}^{(0)} & \mathbf{A}_{22}^{-1}\mathbf{N}_{23}^{(0)} & \mathbf{0} & \mathbf{A}_{22}^{-1}\mathbf{N}_{25}^{(0)}\\
\left(\mathbf{A}_{33}^{(1)}\right)^{-1}\mathbf{N}_{31}^{(0)} & \left(\mathbf{A}_{33}^{(1)}\right)^{-1}\mathbf{N}_{32}^{(0)} & \left(\mathbf{A}_{33}^{(1)}\right)^{-1}\mathbf{N}_{33}^{(1)} & \mathbf{0} & \left(\mathbf{A}_{33}^{(1)}\right)^{-1}\mathbf{N}_{35}^{(1)}\\
\mathbf{0} & \mathbf{0} & \mathbf{0} & \mathbf{A}_{44}^{-1}\mathbf{N}_{44}^{(0)} & \mathbf{A}_{44}^{-1}\mathbf{N}_{45}^{(0)}\\
\left(\mathbf{A}_{55}^{(2)}\right)^{-1}\mathbf{N}_{51}^{(0)} & \left(\mathbf{A}_{55}^{(2)}\right)^{-1}\mathbf{N}_{52}^{(0)} & \left(\mathbf{A}_{55}^{(2)}\right)^{-1}\mathbf{N}_{53}^{(1)} & \left(\mathbf{A}_{55}^{(2)}\right)^{-1}\mathbf{N}_{54}^{(0)} & \left(\mathbf{A}_{55}^{(2)}\right)^{-1}\mathbf{N}_{55}^{(2)}
\end{array}\right].
\]
Finally, Step 4 extracts blocks one level at a time, defining first
\[
\mathbf{P}^{(1)}=\left[\begin{array}{ccccc}
\mathbf{P}_{11}^{(2)} & \mathbf{0} & \mathbf{P}_{13}^{(2)} & \mathbf{0} & \mathbf{P}_{15}^{(2)}\\
\mathbf{0} & \mathbf{P}_{22}^{(2)} & \mathbf{P}_{23}^{(2)} & \mathbf{0} & \mathbf{P}_{25}^{(2)}\\
\mathbf{P}_{31}^{(2)} & \mathbf{P}_{32}^{(2)} & \mathbf{P}_{33}^{(1)} & \mathbf{0} & \mathbf{P}_{35}^{(1)}\\
\mathbf{0} & \mathbf{0} & \mathbf{0} & \mathbf{P}_{44}^{(2)} & \mathbf{P}_{45}^{(2)}\\
\mathbf{P}_{51}^{(2)} & \mathbf{P}_{52}^{(2)} & \mathbf{P}_{53}^{(1)} & \mathbf{P}_{54}^{(2)} & \mathbf{P}_{55}^{(2)}
\end{array}\right]
\]
where the updated blocks are
\begin{alignat*}{1}
\mathbf{P}_{35}^{(1)} & =
\mathbf{P}_{35}^{(2)}-\left(\mathbf{A}_{33}^{(1)}\right)^{-1}
\mathbf{A}_{35}^{(1)}\mathbf{P}_{55}^{(2)}
\\
\mathbf{P}_{53}^{(1)} & =-\left(\mathbf{P}_{35}^{(1)}\right)^{\dagger}
\\
\mathbf{P}_{33}^{(1)} & =\mathbf{P}_{33}^{(2)}-
\left(\mathbf{A}_{33}^{(1)}\right)^{-1}\mathbf{A}_{35}^{(1)}
\mathbf{P}_{53}^{(1)}
.
\end{alignat*}
Note that $\mathbf{P}^{(1)}$ is not skew-Hermitian.
However finalized blocks, such as $\mathbf{P}^{(1)}_{33}$,
$\mathbf{P}^{(1)}_{35}$, $\mathbf{P}^{(1)}_{53}$, and
$\mathbf{P}^{(1)}_{55}$, satisfy the skew-Hermitian
property.
Finally, blocks for regions 1, 2, and 4 are updated, yielding the
matrix $\mathbf{P}^{(0)}$,
\[
\mathbf{P}^{(0)}=\left[\begin{array}{ccccc}
\mathbf{P}_{11}^{(0)} & \mathbf{0} & \mathbf{P}_{13}^{(0)} & \mathbf{0} & \mathbf{P}_{15}^{(0)}\\
\mathbf{0} & \mathbf{P}_{22}^{(0)} & \mathbf{P}_{23}^{(0)} & \mathbf{0} & \mathbf{P}_{25}^{(0)}\\
-\left(\mathbf{P}_{13}^{(0)}\right)^{\dagger} & -\left(\mathbf{P}_{23}^{(0)}\right)^{\dagger} & \mathbf{P}_{33}^{(1)} & \mathbf{0} & \mathbf{P}_{35}^{(1)}\\
\mathbf{0} & \mathbf{0} & \mathbf{0} & \mathbf{P}_{44}^{(0)} & \mathbf{P}_{45}^{(0)}\\
-\left(\mathbf{P}_{15}^{(0)}\right)^{\dagger} & -\left(\mathbf{P}_{25}^{(0)}\right)^{\dagger} & -\left(\mathbf{P}_{35}^{(1)}\right)^{\dagger} & -\left(\mathbf{P}_{45}^{(0)}\right)^{\dagger} & \mathbf{P}_{55}^{(2)}
\end{array}\right].
\]
Blocks for region 4 are satisfying
\begin{alignat*}{1}
\mathbf{P}_{45}^{(0)} & =\mathbf{P}_{45}^{(1)}-\left(\mathbf{A}_{44}^{(0)}\right)^{-1}\mathbf{A}_{45}^{(0)}\mathbf{P}_{55}^{(2)}\\
\mathbf{P}_{44}^{(0)} & =\mathbf{P}_{44}^{(1)}+\left(\mathbf{A}_{44}^{(0)}\right)^{-1}\mathbf{A}_{45}^{(0)}\left(\mathbf{P}_{45}^{(0)}\right)^{\dagger}
\end{alignat*}
Blocks for region 1 are defined by
\begin{alignat*}{1}
\mathbf{P}_{15}^{(0)} & =\mathbf{P}_{15}^{(1)}-\left(\mathbf{A}_{11}^{(0)}\right)^{-1}\mathbf{A}_{15}^{(0)}\mathbf{P}_{55}^{(2)}-\left(\mathbf{A}_{11}^{(0)}\right)^{-1}\mathbf{A}_{13}^{(0)}\mathbf{P}_{35}^{(1)}\\
\mathbf{P}_{13}^{(0)} & =\mathbf{P}_{13}^{(0)}-\left(\mathbf{A}_{11}^{(0)}\right)^{-1}\mathbf{A}_{13}^{(0)}\mathbf{P}_{33}^{(1)}
+
\left(\mathbf{A}_{11}^{(0)}\right)^{-1}\mathbf{A}_{15}^{(0)}
\left(\mathbf{P}_{35}^{(1)}\right)^{\dagger}\\
\mathbf{P}_{11}^{(0)} & =\mathbf{P}_{11}^{(0)}+\left(\mathbf{A}_{11}^{(0)}\right)^{-1}\mathbf{A}_{13}^{(0)}\left(\mathbf{P}_{13}^{(0)}\right)^{\dagger}+\left(\mathbf{A}_{11}^{(0)}\right)^{-1}\mathbf{A}_{15}^{(0)}\left(\mathbf{P}_{15}^{(0)}\right)^{\dagger}
\end{alignat*}
and blocks for region 2
\begin{alignat*}{1}
\mathbf{P}_{25}^{(0)} & =\mathbf{P}_{25}^{(1)}-\left(\mathbf{A}_{22}^{(0)}\right)^{-1}\mathbf{A}_{25}^{(0)}\mathbf{P}_{55}^{(2)}-\left(\mathbf{A}_{22}^{(0)}\right)^{-1}\mathbf{A}_{23}^{(0)}\mathbf{P}_{35}^{(1)}\\
\mathbf{P}_{23}^{(0)} & =\mathbf{P}_{23}^{(0)}
-\left(\mathbf{A}_{22}^{(0)}\right)^{-1}\mathbf{A}_{23}^{(0)}
\mathbf{P}_{33}^{(1)}
+ \left(\mathbf{A}_{22}^{(0)}\right)^{-1}\mathbf{A}_{25}^{(0)}
\left(\mathbf{P}_{35}^{(1)}\right)^{\dagger}
\\
\mathbf{P}_{22}^{(0)} & =\mathbf{P}_{22}^{(0)}+\left(\mathbf{A}_{22}^{(0)}\right)^{-1}\mathbf{A}_{23}^{(0)}\left(\mathbf{P}_{23}^{(0)}\right)^{\dagger}+\left(\mathbf{A}_{22}^{(0)}\right)^{-1}\mathbf{A}_{25}^{(0)}\left(\mathbf{P}_{25}^{(0)}\right)^{\dagger}
\end{alignat*}
All the entries in $\mathbf{P}^{(0)}$ are equal to their corresponding
entries in $\mathbf{G}^{<}$.
\end{document}